\documentclass[12pt,reqno]{amsart}
\usepackage{amsmath,amssymb}
\usepackage[dvips]{graphicx,color}
\usepackage{latexsym,amssymb}
\usepackage{mathrsfs}
\usepackage[abbrev]{amsrefs}

\newtheorem{thm}{Theorem}[section]
\newtheorem{cor}[thm]{Corollaly}
\newtheorem{prop}[thm]{Proposition}
\newtheorem{lem}[thm]{Lemma}

\DeclareMathOperator{\sign}{sign}  
\DeclareMathOperator*{\supp}{supp}

 \newenvironment{pf}
    {{\noindent \bf Proof. }}{\hfill $\Box$}

\numberwithin{equation}{section}
\numberwithin{thm}{section}

\begin{document}

\begin{center}\large \bf 
An application of spectral localization to critical SQG on a ball 
\end{center}

\footnote[0]
{


E-mail: t-iwabuchi@tohoku.ac.jp

}
\vskip5mm

\begin{center}
Tsukasa Iwabuchi 

\vskip2mm

Mathematical Institute, 
Tohoku University\\
Sendai 980-8578 Japan

\end{center}

\vskip5mm

\begin{center}
\begin{minipage}{135mm}
\footnotesize
{\sc Abstract. } 
We study the Cauchy problem for the quasi-geostrophic equations 
in a unit ball of the two dimensional space 
with the homogeneous Dirichlet boundary condition. 
We show the existence, the uniqueness of the strong solution 
in the framework of Besov spaces. 
We establish a spectral localization technique and commutator estimates. 

\end{minipage}
\end{center}

\section{Introduction}

We consider the surface quasi-geostrophic equation in a unit ball. 
\begin{equation}\label{QG1}
\displaystyle 
 \partial_t \theta 
  + (u \cdot \nabla ) \theta + \Lambda _D \theta  =0, 
  \quad u= \nabla ^{\perp} \Lambda _D ^{-1} \theta, 
 \qquad  t > 0 , x \in B, 
\end{equation}
\begin{equation}\label{QG2}
 \theta(0,x) = \theta_0(x) , 
\qquad  x \in B, 
\end{equation}
where 
$B := \{ (x_1, x_2 ) \in \mathbb R^2 \, | \, x_1 ^2 + x_2 ^2 < 1 \}$, 
$\nabla^\perp := (-\partial_{x_2} , \partial_{x_1})$, 
$\Lambda_D$ is the square root of the Dirichlet Laplacian. 
The equations are known as an important model in geophysical fluid dynamics, 
which is derived from general quasi-geostrophic equations in the special case of 
constant potential vorticity and buoyancy frequency (see~\cite{La_1959,Pe_1979}). 
The purpose of this paper is to establish the well-posedness.

Let us recall several known results, where the space is the whole space $\mathbb R^2$. 
If we consider the fractional Laplacian of the order $\alpha$, 
$(-\partial_x^2) ^{\alpha/2}$, with $0 < \alpha \leq 2$, instead of 
$(-\partial_x^2)^{1/2}$, 
then 
the case when $\alpha < 1, \alpha = 1, \alpha > 1$ are called 
sub-critical case, critical case, super-critical case, respectively. 
It is known that the global-in-time regularity is obtained for the 
sub-critical case and the critical case. 
The sub-critical case can be treated, by $L^\infty$-maximum principle, 
and the critical case is delicate. In the critical case, 
the regularity with small data was proved by 
Constantin, Cordoba and Wu~\cite{CoCoWu-2001} 
(see also Constantin and Wu~\cite{ConWu-1999}). 
We also refer on approach in the framework of Besov spaces to 
\cite{Iw-2015,Iw-2020,WZ-2011}
The problem for large data case was solved by Caffarelli and Vasseur~\cite{CaVa-2010}, 
Kiselev, Nazarov and Volberg~\cite{KNV-2007}. 
As another approach, Constantin and Vicol~\cite{CV-2012} proved the 
global regularity by 
nonlinear maximum principles 
in the form of a nonlinear lower bound on the fractional Laplacian. 
On the other hand, in the super-critical case, 
blow-up for smooth solutions is an open problem, 
and the regularity only for small data is known 
(see e.g.~\cite{CotVic-2016}).

In bounded domains with smooth boundary, the equations was introduced by 
Constantin and Ignatova~(\cite{CoIg-2016,CoIg-2017}). Let us focus on the 
critical case. 
Local existence was shown by Constantin and Nguyen~\cite{CoNg-2018-2}, 
and global existence of weak solutions was proved 
by Constantin and Ignatova~\cite{CoIg-2017} for the critical case 
(see also the paper by Constantin and Nguyen~\cite{CoNg-2018} for the inviscid case). 
An interesting question here is how to understand the behavior of the solutions; 
A priori bounds of smooth solutions was obtained by Constantin and 
Ignatova~\cite{CoIg-2016}, 
and interior Lipschitz continuity of weak solutions was studied by Ignatova~\cite{Ig-2019}. 
Recently, Constantin and Ignatova~\cite{CoIg-2020} considered 
the quotient of the solution by the first eigen function to 
investigate near the boundary, and gave a condition to obtain the global regularity 
up to the boundary. 
Stokols and Vasseur~\cite{StVa-2020} constructed global-in-time weak solutions 
with H\"older regularity up to the boundary. 
We should note from the viewpoint of smooth solutions 
that regularity holds for a short time to the best of our knowlegde. 
As for the half space case, the odd reflection reduces the problem to the 
whole spaces case completely, and the analyticity up to the boundary 
is obtained (see~\cite{Iw:preprint2}).

In this paper, we study the existence of solutions for initial data in the 
critical Besov spaces $\dot B^0_{\infty,q}$ associated with the Dirichlet Laplacian, 
where critical space comes from the scaling invariant property in the case when 
the domain is the whole space. Namely, 
the transformation $\theta_\lambda (t,x) = \theta (\lambda t , \lambda x)$ $(\lambda > 0)$ 
maintains the equation \eqref{QG1} and we have 
\[
\| \theta_\lambda (0) \|_X \simeq \| \theta(0) \|_X 
\quad \text{ for all } \lambda > 0,
\]
for $X = L^\infty(\mathbb R^2), \dot H^\frac{2}{p}_p(\mathbb R^2), 
\dot B^\frac{2}{p}_{p,q}(\mathbb R^2)$. It would be natural that these spaces on domains  
have some critical structure locally in time at least. 
We prove the existence of local solutions for arbitrary data and 
global solutions for small data.

We state our main result for initial data in $\dot B^0_{\infty,1}$ to explain 
the essence of this paper simply, and mention that 
$\dot B^0_{\infty,q}$, $q > 1$, case follows as well as the whole space case. 
We also see that in the case when $q = 1$ the functions in $\dot B^0_{\infty,1}$ 
is continuous up to the boundary, and the boundary condition is understood by 
the boundary value of continuous functions. 

We introduce the definition of Besov spaces. 
Let $\phi_0 $ be such that $\phi_0 \in C_0^\infty (\mathbb R)$ and 
\[
{\rm supp \, } \phi_0 \subset [2^{-1} , 2], \quad 
\phi_0 (\lambda) = \phi_0 \left( \dfrac{\lambda}{2^j} \right) , 
\quad 
\sum _{j \in \mathbb Z} \phi_0 \left( \dfrac{\lambda}{2^j} \right)   = 1  
\text{ for any } \lambda > 0 ,
\]
and we define 
\[
\phi_j(\lambda) = \phi_0 \left( \dfrac{\lambda}{2^j} \right) , \quad \lambda \in \mathbb R. 
\]
For $s \in \mathbb R$ and $1 \leq p,q \leq \infty$, $\dot B^s_{p,q} = \dot B^s_{p,q}(\Lambda _D)$ is defined by 
\[
\dot B^s_{p,q} = \dot B^s_{p,q}(\Lambda _D)
= \left\{ f \in \mathcal Z_D'   \, \Big| \, 
 \| f \|_{\dot B^s_{p,q}} = \Big\{ \sum_{j \in \mathbb Z} \Big( 2^{sj} \| \phi_j(\Lambda _D) f \|_{L^p} \Big) q
 \Big\}^{\frac{1}{q}} < \infty \right\},
\]
where $\mathcal Z'_D$ is a space of tempered distribution, 
and we explain the precise definitions in subsection~\ref{subsec:besov}.

\begin{thm}
\label{thm:1} 
For every 
$\theta _0 \in  \dot B^0_{\infty,1} (\Lambda _D)$, there exists $T > 0$ such that 
the equation \eqref{QG1} with the initial condition \eqref{QG2} has a unique solution 
$\theta$ such that 
\begin{gather*}
\theta \in C([0,T ], \dot B^0_{\infty,1}) \cap L^1(0,T ; \dot B^1_{\infty,1}), 
\quad 
\partial _t \theta \in L^1 (0,T ; \dot B^1_{\infty,1}),
\end{gather*}
and $\theta$ is continuous with respect to $t \geq 0 , x \in B$, 
and 
\[
\lim_{|x| \to 1} \theta (t,x) = 0 \quad \text{for } t \geq 0 .
\]
If $\theta_0$ is small in $\dot B^0_{\infty,1}$, then the solution $\theta$ exists globally in time. 
\end{thm}

Let us give some remarks to prove theorems. 
Our idea is to establish a method based on that in the case of the whole space (\cite{WZ-2011,Iw-2020}), replacing Fourier transformation with spectral decomposition,  
and the difference is that we have several problems,  how to control or understand 
boundary value of functions satisfying the Dirichlet boundary condition. 
Our starting point is the spectral multiplier theorem, which is boundedness of the 
operator $\varphi (-\Delta _D)$ in $L^p$, $1 \leq p \leq \infty$ for all 
$\varphi$ in the Schwartz class in the real line, and 
our ingredient is applying spectral localization (see Proposition~\ref{prop:0906-3}) 
and commutator estimates (see Proposition~\ref{prop:0906-2}). 
The spectral localization in this paper is a bound from below, more precisely, 
at a maximum point $x_0$ of $|\psi_j(\Lambda_D) f| $
\[
\Lambda _D \psi_j(\Lambda_D) f(x_0)  \sign f(x_0)\geq  c 2^j \| \psi_j (\Lambda _D)f \|_{L^\infty},
\]
where $\psi _j(\Lambda)$ is an operator restricting the spectrum around $2^j$. 
We remark that this is possible for other domains. 
This kind of localization when the domain is the whole space is established in \cite{WZ-2011}, and it can be generalized in domains as in Proposition~\ref{prop:0906-3}. 
The commutator estimate in this paper is a bilinear estimate, 
\[
\sum_{j \in \mathbb Z} 
\Big\| \Big[\nabla ^\perp \Lambda _D^{-1} f \cdot \nabla , \,\, \psi_j(\Lambda _D) 
  \Big] g
\Big\|_{L^\infty}
\leq C \| f \|_{\dot B^1_{\infty,1}} \| g \|_{\dot B^0_{\infty,1}} .
\]
However it seems very difficult to estabilish by Littlewood Paley dyadic decomposition as in 
the whole space, while we can avoid by using a resolution of unity such that 
\[
1 = \sum_{j \in \mathbb Z} \left( \frac{1}{1+2^{-2j-2}\lambda ^2} - \frac{1}{1+2^{-2j} \lambda ^2} \right)
=: \sum _{j \in \mathbb Z}\psi_j(\lambda), \quad 
\lambda > 0 ,
\]
since the resolvent have a property that 
\[
\begin{split}
& 
\Big[\nabla ^\perp \Lambda _D^{-1} f \cdot \nabla , \,\, \frac{1}{1-2^{-2j} \Delta _D} 
  \Big] g
\\
  =
&  \frac{-2^{-2j}}{1-2^{-2j}\Delta _D} 
  \Big( (\nabla ^{\perp}\Delta _D \Lambda _D^{-1}  f \cdot \nabla) 
      + 2(\nabla ^\perp\nabla \Lambda _D^{-1} f \cdot \nabla ) \nabla \Big)
   \dfrac{1}{1-2^{-2j}\Delta _D} g, 
\end{split}
\]
and the right hand side is justfied if $g$ does not have weak derivatives, 
where $-\Delta _D$ denotes the Dirichlet Laplacian.  
We also see the equivalency of the norm defined by the dyadic decomposition and the resolvent 
for $s$ close to zero, and the commutator estimate can be expected. 
It is also important in the argument above 
that if $f,g $ satisfy the Dirichlet boundary condition and are smooth, then 
$( \nabla^\perp\Lambda_D^{-1}f \cdot \nabla )g$ does and 
\[
( \nabla^\perp\Lambda_D^{-1}f \cdot \nabla )g 
= (1-2^{-2j}\Delta_D)^{-1}  (1-2^{-2j} \Delta _D)
\Big( ( \nabla^\perp\Lambda_D^{-1}f \cdot \nabla )g \Big) .
\]
Finally, we mention basic tools; maximum regularity estimates in 
Lemma~\ref{lem:0906-1} below~(see also~\cite{Iw-2018}), 
bilinear estimates in Proposition~\ref{prop:0830-10}. 
In this paper, 
we give a simple proof for boundedness of second derivative near $L^\infty$ space in Lemma~\ref{lem:0830-6} by using the explicit formula of 
$(-\Delta _D)^{-1}$ in a ball, but it would be possible to consider 
other smooth bounded domains.

\bigskip

It is also possible to obtain a similar result with 
initial data in $\dot B^0_{\infty,q}$ with $q > 1$, by introducing spaces, whose norms are defined by 
\[
\|f\| _{\widetilde L^p (0,T ;\dot B^s_{\infty,q})} := 
 \Big\|  
 \Big\{ 2^{sj} \| \phi_j(\Lambda _D) f \|_{L^1 (0,\infty;L^p(\mathbb R^2_+))} 
 \Big\}_{j \in \mathbb Z}
\Big\| _{\ell^q(\mathbb Z)} .
\]
\begin{thm}
\label{thm:2} 
\begin{enumerate}
\item[\rm (i)]
Let $1 \leq q < \infty$. 
For every 
$\theta _0 \in  \dot B^0_{\infty,q} (\Lambda _D)$, there exists $T > 0$ such that 
the equation \eqref{QG1} with the initial condition \eqref{QG2} has a unique solution 
$\theta$ such that 
\begin{gather*}
\theta \in C([0,T ], \dot B^0_{\infty,q}) 
\cap \widetilde L^\infty (0,T ; \dot B^0_{\infty,q})
\cap \widetilde L^1(0,T ; \dot B^1_{\infty,q}), 
\quad 
\partial _t \theta \in \widetilde L^1 (0,T ; \dot B^1_{\infty,q}),
\end{gather*}
If $\theta_0$ is small in $\dot B^0_{\infty,q}$, then the solution exists globally in time. 

\item[\rm (ii)]
Let $ q = \infty$. For every $\theta _0 \in \dot B^0_{\infty, \infty}$ such that 
$\| \phi_j(\Lambda _D) \theta_0 \|_{L^\infty} \to 0$ $(j \to \infty)$, 
the same existence result holds. 
\end{enumerate}

\end{thm}

\bigskip 

This kind of local existence  when the domain is the whole space is established 
by Wang-Zhang~\cite{WZ-2011} (see also \cite{Iw-2020}). 
It is possible to modify the proof of Theorem~\ref{thm:1} to handle the case 
when $q > 1$ in a similar way to \cites{Iw-2020,WZ-2011}. 
We left the proof for readers. 

\bigskip 

This paper is organized as follows. 
In section 2, we recall the definition of Besov spaces associated with the Dirichlet 
Laplacian, and several properties for the boundary value of functions, 
such as spectral localization, commutator estimates. 
In section 3, we prove Theorem~\ref{thm:1}. 
In appendix, we discuss the equivalence of the Besov norms 
defined by the dyadic decomposition and the resolvent when the regularity 
is close to $0$.

\bigskip

\noindent 
{\bf Notations.} We denote by $-\Delta _D$ the Dirichlet Laplacian on $L^2 (\Omega)$ 
and on distribution space $\mathcal Z_D'$ defined in section 2. 
We write $x=(x_1,x_2)$. 
Let $\{ \phi_j \}_{j \in \mathbb Z}$ be the dyadic decomposition of the 
unity such that $\phi_j$ is a non-negative function in $C_0 ^\infty (\mathbb R)$ and 
\[
{\rm supp \, } \phi_0 \subset [2^{-1} , 2], \quad 
\phi_j (\lambda) = \phi_0 \left( \dfrac{\lambda}{2^j} \right) , 
\quad 
\sum _{j \in \mathbb Z} \phi_j (\lambda)  = 1  
\text{ for any } \lambda > 0 . 
\]
$\{ \psi_j \}_{j \in \mathbb Z}$ is another resolution of identity such that 
\[
1 = \sum_{j \in \mathbb Z} \left( \frac{1}{1+2^{-2j-2}\lambda ^2} - \frac{1}{1+2^{-2j} \lambda ^2} \right)
=: \sum _{j \in \mathbb Z}\psi_j(\lambda), \quad 
\lambda > 0 . 
\]
We use the following notations for norms of spaces in space and time as follows. 
\[
\begin{split}
\| f \|_{\dot B^s_{p,q}(\Lambda_D)} = 
&
\Big\|  
 \Big\{ 2^{sj} \| \phi_j(\Lambda _D) f \|_{L^p} 
 \Big\}_{j \in \mathbb Z}
\Big\| _{\ell^q(\mathbb Z)},
\\
\| f \|_{L^r(0,\infty; X)} 
= 
&\big\| \| f(t) \|_{X} \big\|_{L^r(0,\infty)} , 
\qquad X = L^p (\mathbb R^2_+), \dot B^s_{p,q} (\Lambda _D). 
\end{split}
\]
We define the Sobolev spaces $H^s$ by the Besov spaces
\[
H^s = \dot B^s_{2,2}, \quad s \in \mathbb R.
\]
We write the domain of the functions only 
when the function space on the whole space is 
used, for instance $ B^s_{p,q}(\mathbb R^2)$, where we will use the theory on $\mathbb R^2$. 
When the domain is the ball $B$, then we omit it. 
The kernel of $(-\Delta _D)^{-1}$ (see e.g., \cite{Evans}) is defined by 
\[
(-\Delta _D)^{-1} (x,y) 
= \frac{1}{2\pi}\log |x-y| - \Phi(x,y), 
\]
where 
\begin{equation}\label{0906-10}
\Phi (x,y) = \dfrac{1}{2\pi} \log  \Big( |x| \Big|y - \dfrac{x}{|x|^2} \Big| \Big) . 
\end{equation}
We write 
\[
B:=\{ x \in \mathbb R^2 \, | \, |x|<1 \}, \quad 
B^c := \{ x\in \mathbb R^2  \, | \, |x| \geq 1\},
\]
and $\chi_B$, $\chi_{B^c}$ the characteristic functions on 
$B, B^c$, respectively. 
On the whole space $\mathbb R^2$, we denote by 
$-\Delta_{\mathbb R^2} , \Lambda_{\mathbb R^2}$, 
the Laplacian and the square root of the Laplacian defined by the Fourier transform 
defined in the space of tempered distribution.  
We denote by $-\Delta_{D}, \Lambda_D$, 
the Dirichlet Laplacian and the square root of the Dirichlet Laplacian.

\section{Preliminary}

In subsection 2.1, we recall the definition of Besov spaces in \cite{IMT-2019}. 
In subsection 2.2, we introduce spectral multiplier theorem together with derivative 
estimates and smoothing property such as maximum regularity for $e^{-t \Lambda _D}$. 
In subsection 2.3, inequalities for the spectral localization is estabilished. 
In subsection 2.4, the commutator estimates are proved.

\subsection{Besov spaces}\label{subsec:besov}

We recall the definition of the Besov spaces (see~\cite{IMT-2019}). 
We start by defining the Dirichlet Laplacian $-\Delta _D$, and 
spaces of test functions, $\mathcal Z$, of homogeneous type. 
We here notice that the infimum of the spectrum is strictly positive, 
since we consider the bounded domain with the Dirichlet condition, 
and the spaces of homogeneous and non-homogeneous types are equivalent. 
We just adopt the homogeneous type for a simple notation in our proof.

\bigskip

\noindent 
{\bf Definition}. 
(i)
Let $-\Delta _D $ be the Dirichlet Laplacian on $L^2 (B)$ defined by 
\[
\begin{cases}
D(-\Delta _D) := \{ f \in H^1_0(B)  \, | \, 
 \Delta f \in L^2 (B\}, 
\\
-\Delta _D f := -\Delta f =
- \left( \dfrac{\partial^2}{\partial x_1 ^2} f 
  + \dfrac{\partial^2}{\partial x_2 ^2} f 
  \right) , 
\quad f \in D(-\Delta _D). 
\end{cases}
\]

\noindent 
(ii) Let $\mathcal Z _D $ 
be a space of test functions  such that 
\[
\mathcal Z_D := \{ f \in L^1 \cap L^2 \, | \, 
   q_m (f) < \infty \text{ for all } m \in \mathbb N  \}, 
\]
where 
\[
q_m(f):= 
\sup _{ j  \in \mathbb Z} 2^{m |j|} \| \phi_j (\Lambda _D) f \|_{L^1} . 
\]

\noindent 
(iii) Let $\mathcal Z_D '$ be the topological duals of $\mathcal Z_D$.

\bigskip 

It was proved in \cite{IMT-2019} that 
the space $\mathcal Z_D$ is a Fr\'echet space, 
and can regard their duals $ \mathcal Z_D'$  
as distribution spaces, 
which are variants of the space of the tempered distributions and the quatient space 
by the polynomials in the whole space. 
We define Besov spaces associated with the Dirichlet Laplacian 
on the unit ball as follows. 

\bigskip

\noindent 
{\bf Definition}. 
Let $s \in \mathbb R$ and $1 \leq p,q \leq \infty$. 
$\dot B^s_{p,q} = \dot B^s_{p,q}(\Lambda_D)$ is defined by 
\[
\dot B^s_{p,q} = \dot B^s_{p,q}(\Lambda_D) :=  
\{ f \in \mathcal Z_D' \, | \, 
   \| f \|_{\dot B^s_{p,q}(\Lambda_D) } < \infty\}, 
\]
where 
\[
\| f \|_{\dot B^s_{p,q}(\Lambda_D)} 
:= 
\Big\|  
 \Big\{ 2^{sj} \| \phi_j(\Lambda _D) f \|_{L^p(\mathbb R^2_+)} 
 \Big\}_{j \in \mathbb Z}
\Big\| _{\ell^q(\mathbb Z)}.
\]

\bigskip 

It is proved that $\dot B^s_{p,q}(\Lambda _D) $ 
is a Banach space and satisfies standard properties such as  lift properties, 
embedding theorems of Sobolev type 
as well as the whole space case. 
We here recall the uniform boundedness of the frequency restriction operator 
$\phi_j(\Lambda_D)$ and some fundamental property of the Besov spaces 
for our purpose of this paper. This is possible, since 
operators $\phi_j(\Lambda)$ ($j \in \mathbb Z$) restrcting the spectrum 
are uniformly bounded in $L^p$ for all $1 \leq p \leq \infty$ 
(see~Lemma~\ref{lem:0906-7} for more details).

\bigskip 

We here write several properties which are needed in the proof.

\begin{lem}\label{lem:0906-8}
Let $|s| < 2$ and $1 \leq p,q \leq \infty$. Then 
\[
f = \sum _{ j \in \mathbb Z} \phi_j(\Lambda _D) f 
\quad \text{ in } \mathcal Z'_D ,
\qquad  
\| f \|_{\dot B^s_{p,q}} 
\simeq 
\left\{ \sum _{ j \in \mathbb Z} 
 \Big( 2^{sj} \| \psi_j (\Lambda _D) f \|_{L^p} \Big) ^q
\right\}^{\frac{1}{q}}, 
\]
for all $f \in \dot B^s_{p,q}$. 
\end{lem}

We prove Lemma~\ref{lem:0906-8} in Appendix~\ref{Appen_1}.

\begin{lem}Let $1 \leq p,q,r \leq \infty$. Then 
\begin{equation}
\label{0902-6}
\| f \|_{L^p} \leq \| f \|_{ \dot B^0_{p,1}}, \|\nabla f\|_{L^\infty} \leq C \| f\|_{ \dot B^1_{p,1}},
\end{equation}
\begin{equation}\label{0915-1}
\| \Lambda_D ^s f \|_{\dot B^0_{p,q}} 
\leq C \| f \| _{\dot B^s_{p,q}} \quad \text{for } s \in \mathbb R , 
\end{equation}
\begin{equation}\label{0902-8}
\| f \|_{\dot B^0_{p,q}} 
\leq C \| f \|_{\dot B^{2(\frac{1}{r}-\frac{1}{p})}_{r,q}},
\end{equation}
\begin{equation}
\label{0902-7}
\| f \|_{\dot  B^0_{\infty,1}} \leq C \| f \|_{\dot B^s_{\infty,q}} \text{ for } s > 0 .
\end{equation}
\end{lem}

\begin{pf}
The first inequality of \eqref{0902-6} is obtained by the resolution of the identity in Lemma~\ref{lem:0906-8} 
and the triangle inequality, and 
the second inequality is proved by the resolution and the gradient estimate \eqref{0830-2}. 
The lifting property, the embedding theorem is already known in \cite{IMT-2019}. 
The validity of the last inequality \eqref{0902-7} is due to the infimum of the spectrum being positive 
and an elementary boundedness in the sequence spaces. 
\end{pf}

\begin{lem}\label{lem:0906-4}
\begin{enumerate}
\item[\rm (i)] Every $f \in \dot B^0_{\infty,1}$ is regarded as a continuous function 
up to the boundary and $f \equiv 0$ on the boundary. 
\item[\rm (ii)] 
Let $f ,g \in L^\infty$ and $f_k = \phi_k(\Lambda _D) f$, 
$g_l := \phi_l (\Lambda _D) g$ for $k,l \in \mathbb Z$. 
Then $(\nabla ^{\perp}\Lambda _D f_k \cdot \nabla ) g_l$ is regarded as 
a continuous function up to the boundary and is equal to zero on the boundary. 
\end{enumerate}
\end{lem}

We give a direct proof of Lemma~\ref{lem:0906-4} by using the formula of $(-\Delta _D)^{-1}$ 
in appendix~\ref{appen_2}, since the proof seems elementary. 
One can also find the orthogonality due to $\nabla ^\perp$ and $\nabla$ for 
for functions in $H^1_0$ on smooth bounded domain in \cite{CoNg-2018-2}.

\subsection{Spectral multiliers and smoothing property of $e^{-t\Lambda _D}$}

We recall boundedness of the spectral multipliers and gradient estimates. 
We mainly refer \cite{IMT-2018}, but there are a plenty of literature on this field, 
and one can refer to \cite{Ou_2005,IMT-2021,ThOuSi-2002} for the theory. 

\begin{lem} {\rm (\cite{IMT-2018})} {\rm (Boundedness of spectral multiplers)}  
\label{lem:0906-7} 
Suppose that $\psi $ and its all derivatives are bounded and that 
$\varphi$ belongs to the Schwartz class in the real line and $1 \leq p \leq \infty$. 
Then $\delta , C > 0$ exist such that 
\begin{equation}\label{0902-1}
\| \varphi (2^{-2j} (-\Delta _D) \|_{L^p \to L^p} 
\leq C \| (1-\Delta_{\mathbb R})^{\frac{d+1}{4}+ \delta}
           \varphi(\cdot) \|_{L^2 (\mathbb R)} .
\end{equation}
Furthermore, 
\begin{equation}\label{0901-1}
\| \psi(-\Delta _D) \varphi (2^{-2j} (-\Delta _D) \|_{L^p \to L^p} 
\leq C \| (1+|\cdot|^2)^{\frac{3d}{8}+\frac{d}{4}+ \delta } (1-\Delta_{\mathbb R})^{\frac{d+1}{4}+ \delta}
          \psi(2^{2j} \cdot) \varphi(\cdot) \|_{L^2 (\mathbb R)} .
\end{equation}
\end{lem}

We recall maximum regularity estimate. 

\begin{lem}{\rm (\cite{Iw-2018})}\label{lem:0906-1}
For every $\theta_0 \in \dot B^0_{\infty,1}$ 
\begin{equation}\label{0902-2}
\| e^{-t\Lambda}\theta_0 \|_{L^\infty (0,T ; \dot B^0_\infty,1) \cap L^1 (0,T ; \dot B^0_{\infty,1})} 
\leq C \| \theta_0 \|_{\dot B^0_{\infty,1}} . 
\end{equation}
If $u \in C([0,T] , \dot B^0_{\infty,1}) \cap L^1 (0,T ;\dot B^1_{\infty,1})$ 
and $f \in L^1 (0,T ; \dot B^0_{\infty,1})$ satisfy 
$\partial _t u \in L^1 (0,T ; \dot B^1_{\infty,1})$, 
$\partial _t u + \Lambda u = f$, then 
\begin{equation}\label{0902-3}
\| u \|_{L^\infty(0,T;\dot B^0_{\infty,1}) \cap L^1 (0,T ; \dot B^1_{\infty,1})} 
\leq C \|u(0) \|_{\dot B^0_{\infty,1}} + C \| f \|_{L^1 (0,T ; \dot B^0_{\infty,1} )} .
\end{equation}
\end{lem}

We also use the boundedness of the resolvent. 

\begin{lem}For every $1 \leq p \leq \infty$
\begin{equation}\label{0902-12}
\sup _{j \in \mathbb Z} \left\| (1- 2^{-2j} \Delta _D)^{-1} \right\|_{L^p \to L^p} < \infty.
\end{equation}
\end{lem}

\begin{pf}
Let $f \in L^p $ and $\widetilde \phi_0 \in C_0^\infty (\mathbb R)$ be such that 
\[
\supp \widetilde \phi_0 \subset [2^{-1}, 2], \quad 
\widetilde \phi_0 (\lambda^2) = \sum _{j \leq 0} \phi_j(\lambda), \quad \lambda > 0.
\]
We use the resolution 
\[
1 = \widetilde \phi_0 (2^{-2j}\lambda ^2) + \sum _{k > j} \phi_j (\lambda) , 
\quad \text{ for } \lambda > 0 .
\]
By $(1+\lambda^2)^{-1} \widetilde \phi_0 \in C_0^\infty (\mathbb R)$ 
and the boundedness of the spectral multipliers \eqref{0902-1}, 
\[
\begin{split}
& 
\| (1-2^{-2j} \Delta _D)^{-1} f \|_{L^p} 
\\
\leq 
& 
\| (1-2^{-2j} \Delta _D)^{-1} \widetilde \phi_0( - 2^{-2j}\Delta _D) f \|_{L^p} 
 + \sum _{ k>j} \| (1-2^{-2j} \Delta _D)^{-1} \phi_k( \Lambda _D )f \|_{L^p} 
\\
\leq 
& 
C\|f \|_{L^p} 
+ \sum _{ k>j}\frac{C}{1+2^{-2j+2k}} \| f \|_{L^p} 
\leq C \| f \|_{L^p} .
\end{split}
\]
\end{pf}

We use the boundedness of the derivatives. 
\begin{lem}{\rm (\cite{IMT-2018})} \label{lem:0830-5} 
Let $m = 0,1,2, \cdots, $ and $1 \leq p \leq \infty$. Then 
\begin{gather}\label{0830-3}
\| \Lambda_D ^m \phi_j(\Lambda _D) f \|_{L^p} 
\leq C 2^{mj} \|  \phi_j(\Lambda _D) f \|_{L^p}, 
\\ 
\label{0830-2}
\| \nabla \phi_j(\Lambda _D)f \|_{L^p} +
\|  \phi_j(\Lambda _D) \nabla f \|_{L^p} \leq C 2^j \| \phi_j(\Lambda _D) f\|_{L^p}.
\end{gather}
\end{lem}

\begin{lem}\label{lem:0830-6}
{\rm (i)} There exists a constant $C > 0$ such that 
\begin{equation}\label{0819-1}
\| \nabla ^2 (-\Delta_D )^{-1} f \|_{L^\infty} 
\leq C \| f \|_{\dot B^0_{\infty,1}} , 
\quad \text{ for } f \in \dot B^0_{\infty,1} . 
\end{equation}
When $1 < p < \infty$, 
\begin{equation}\label{0819-1-2}
\| \nabla ^2 (-\Delta_D )^{-1} f \|_{L^p} 
\leq C \| f \|_{L^p} , 
\quad \text{ for } f \in L^p . 
\end{equation}

\bigskip 

\noindent 
{\rm (ii)} 
Let $f,g \in L^\infty$ and $f_k := \phi_k (\Lambda_D)f$, $g_l := \phi_l (\Lambda _D) g$ $(k,l \in \mathbb Z)$. 
Then 
\begin{equation}\label{0830-4}
\nabla ^\perp f_k \cdot \nabla g_l \in H^1_0 \quad \text{ and } \quad 
(-\Delta _D) \Big( \nabla ^\perp f_k \cdot \nabla g_l  \Big) \in L ^2 . 
\end{equation}
\end{lem}

\begin{pf}
By $\Phi(x,y) = \Phi(y,x)$ (see \eqref{0906-10} for the definition) 
and a change of variable $y \mapsto y/|y|^2$, we write 
\[
\begin{split}
& (-\Delta_D)^{-1} f (x)
\\
=& \frac{1}{2\pi} \int_{\{|x|<1 \}} \Big( \log|x-y| \Big) f(y) ~dy 
-\frac{1}{2\pi} \int_{B^c} \Big( \log |y|^{-1} + \log|x-y| \Big) f\Big(\dfrac{y}{|y|^2} \Big) \frac{dy}{|y|^4}
\\
=& \frac{1}{2\pi} \int_{\mathbb R^2} 
  \Big( \log|x-y| \Big) \Big( \chi_{B} f(y) -  \dfrac{\chi_{B^c} }{|y|^4}  f\Big(\dfrac{y}{|y|^2} \Big)  \Big) 
  dy
+ \frac{1}{2\pi} \int_{B^c}\dfrac{\log|y|}{|y|^4} f\Big(\dfrac{y}{|y|^2} \Big) dy,
\end{split}
\]
and an extention to $\mathbb R^2$
\[
F(y) := \chi_{B} f(y) -  \dfrac{\chi_{B^c} }{|y|^4}  f\Big(\dfrac{y}{|y|^2} \Big), \quad y \in \mathbb R^2.
\]
By considering $\mathbb R^2$ and the boudedness of the Riesz transform (see e.g. Stein~\cite{Stein_1970}), 
\begin{equation}\label{0830-1}
\begin{split}
\| \nabla^2 (-\Delta_D)^{-1} f \|_{L^\infty} 
\leq 
& 
 \| \nabla ^2 (-\Delta_{\mathbb R^2})^{-1} F \|_{L^\infty (\mathbb R^2)}
\leq \sum_{j \in \mathbb Z} \|\nabla ^2 (-\Delta_{\mathbb R^2})^{-1} \phi_j (\Lambda_{\mathbb R^2}) F\|_{L^\infty(\mathbb R^2)}
\\
\leq 
& C \| F \|_{ \dot B^0_{\infty,1}(\mathbb R^2)} .
\end{split}
\end{equation}
The real interpolation $\dot B^0_{\infty,1}(\mathbb R^2) 
= (BMO^{-2}(\mathbb R^2) , BMO^{2}(\mathbb R^2) )_{\frac{1}{2},1}$ 
(see e.g. \cite{Tri_1983}) implies that 
\[
\begin{split}
\| F \|_{ \dot B^0_{\infty,1}(\mathbb R^2)}
\leq \int_0^\infty t^{-\frac{1}{2}} 
  \inf_{F = F_1+F_2} 
   \Big( \left\| (-\Delta_{\mathbb R^2})^{-1}F_1 \right\|_{BMO(\mathbb R^2)} 
              + t \left\| (-\Delta _{\mathbb R^2}) F_2 \right\|_{BMO(\mathbb R^2)}
   \Big)  \frac{dt}{t} .
\end{split}
\]
Here we restrict the decomposition $F = F_1 + F_2$ to $F_1, F_2$ such that 
\[
F_j = \chi_{B} f_j(y) -  \dfrac{\chi_{B^c} }{|y|^4}  f_j\Big(\dfrac{y}{|y|^2} \Big), \quad j = 1,2,
\]
where $f_1 ,f_2$ are functions satisfying $f = f_1 + f_2$ in the ball $B$. We then have 
\[
\begin{split}
\left\| (-\Delta_{\mathbb R^2})^{-1}F_1 \right\|_{BMO(\mathbb R^2)} 
=&  \left\| (-\Delta_{\mathbb R^2})^{-1}F_1 
        + \frac{1}{2\pi} \int_{B^c}\dfrac{\log|y|}{|y|^4} f_1\Big(\dfrac{y}{|y|^2} \Big) dy
    \right\|_{BMO(\mathbb R^2)}
\\
= & \left\| \widetilde {(-\Delta_D)^{-1} f_1} \right\|_{BMO(\mathbb R^2)}
\leq C \| (-\Delta_D)^{-1}f_1 \|_{L^\infty} 
\leq C \| f_1 \|_{\dot  B^{-2}_{\infty,1}} ,
\end{split}
\]
where 
\[
\widetilde {(-\Delta_D)^{-1} f_1}(x) 
= 
\begin{cases}
\Big((-\Delta_D)^{-1} f_1 \Big)(x) & \text{ for } x \in B , 
\\
- \Big((-\Delta_D)^{-1} f_1 \Big)\Big( \dfrac{x}{|x|^2} \Big)  & \text{ for } x \in B^c , 
\end{cases}
\]
and 
\[
\left\| (-\Delta _{\mathbb R^2}) F_2 \right\|_{BMO(\mathbb R^2)}
\leq C ( \| f_2 \|_{L^\infty} + \| \nabla f_2\|_{L^\infty} + \| (-\Delta_D) f_2 \|_{L^\infty} )
\leq C \| f_2 \|_{ \dot B^2_{\infty,1}} .
\]
By the real interpolation $ \dot B^0_{\infty,1} = ( \dot B^{-2}_{\infty,1} ,  \dot B^2_{\infty,1})_{\frac{1}{2},1}$, 
we conclude  
\[
\begin{split}
\| F \|_{ \dot B^0_{\infty,1}(\mathbb R^2)}
\leq C \int_0^\infty t^{-\frac{1}{2}} 
  \inf_{f = f_1 + f_2} 
   \Big( \left\| f_1 \right\|_{\dot  B^{-2}_{\infty,1}} 
              + t \left\|f_2 \right\|_{ \dot B^2_{\infty,1}}
   \Big)  \frac{dt}{t} 
   \leq C \| f \|_{\dot  B^0_{\infty,1}},
\end{split}
\]
which proves \eqref{0819-1}. 
The boundedness in $L^p$, $1 < p < \infty$, follows from a similar argument and the boundedness of the Riesz transform (see e.g. \cite{Stein_1970}) in the whole space.

Next, we write 
\[
f_k = (-\Delta_D)^{-1} \Big( -\Delta _D f_k \Big) , \quad 
g_l = (-\Delta_D)^{-1} \Big( -\Delta _D g_l \Big),
\]
and the boundary value of $\nabla ^\perp f_k \cdot \nabla g_l$ can be checked by the explicit formula 
of $(-\Delta_D)^{-1} (x,y) $ and its derivatives (or the orthogonality of $\nabla ^\perp$ and $\nabla$ 
(see~\cite{Constantin ??})), and it yields that 
$\nabla ^\perp f_k \cdot \nabla g_l$ is continuous up to the boundary and 
$\nabla ^\perp f_k \cdot \nabla g_l = 0$ on the boundary. 
We can also see from a similar argument to \eqref{0830-1} and the spectral multiplier theorem~\eqref{0830-3}that 
\[
\begin{split}
\| \nabla ^2 f_k \|_{L^4} \leq 
& 
\left\| \nabla ^2 (-\Delta _{\mathbb R^2})^{-1} \Big( \chi_B (-\Delta _D)f_k - \frac{\chi_{B^c} (-\Delta _D)f_k (\frac{\cdot}{|\cdot|^2})}{|y|^4}  \Big)
\right\|_{L^4(\mathbb R^2)}
\\
\leq 
& C \| (-\Delta _D)f_k \|_{L^4} 
\leq C 2^{2k} \| f_k \|_{L^4},
\end{split}
\]
and $g_l$ satisfies the same inequality, with replacing $2^{2k}$ by $2^{2l}$. 
By the H\"older inequality and the gradient estimate~\eqref{0830-2}, we conclude that 
\[
\begin{split}
\| \nabla (\nabla^{\perp}f_k \cdot \nabla g_l) \|_{L^2} 
\leq & \| \nabla ^2 f _k \|_{L^4} \| \nabla g_l \|_{L^4} + \| \nabla f_k \|_{L^4} \| \nabla ^2 g_l \|_{L^4}
\\
\leq & C (2^{2k +l} + 2^{k+2l}) \| f_k \|_{L^4} \| g_l \|_{L^4} < \infty ,
\end{split}
\]
and $\nabla ^\perp f_k \cdot \nabla g_l \in H^1_0 $. Also, $(-\Delta _D) f \in L^2 $ follows from 
\[
\begin{split}
(-\Delta _D) (\nabla^{\perp}f_k \cdot \nabla g_l) 
=& \nabla^{\perp}(-\Delta _D)f_k \cdot \nabla g_l  
 - 2 \nabla \nabla^{\perp}f_k \cdot \nabla \nabla g_l  
 + \nabla^{\perp}f_k \cdot \nabla (-\Delta_D )g_l , 
\\
\| (-\Delta _D) (\nabla^{\perp}f_k \cdot \nabla g_l)  \|_{L^2} 
\leq & C (2^{3k +l} + 2^{2k + 2l}+ 2^{k+3l}) \| f_k \|_{L^4} \| g_l \|_{L^4} < \infty. 
\end{split}
\]
\end{pf}

\subsection{Spectral localization}

\begin{prop}\label{prop:0906-3}
Let $\psi_j (\lambda) := (1+2^{-2j-2}\lambda ^2)^{-1} - (1+2^{-2j}\lambda ^2)$. 
There exists $c > 0$ such that for every $f \in L^\infty$ and $j \in \mathbb Z$, we have 
at a maximum point $x_0$ of $ |\phi_j(\Lambda _D) f| $  
\begin{equation}\label{0827-1}
\Big| \Big( \Lambda_D \psi_j (\lambda _D) f \Big)(x_0) \Big|  
\geq c 2^j  \| \psi_j (\Lambda _D) f \|_{L^\infty} . 
\end{equation}

\end{prop}

\begin{pf}
We write $f_j =  \psi_j(\Lambda _D) f$ for the sake of simplicity, and 
we may assume positivity $f_j(x_0) = \| f_j \|_{L^\infty} \geq 0$ 
at the maximum point $x_0$, 
unless we may consider the case when $f_j(x_0) $ is negative and it suffices to replace $f$ by $-f$. 
We choose a constant $c_1>0$ such that 
\[
\lambda = c_1 \int_0^\infty t^{-\frac{3}{2}} (1-e^{-t\lambda ^2}) dt, 
\quad \text{ for } \lambda > 0 ,
\]
and write 
\[
\Lambda _D f_j (x_0) = c_1 \int_0^\infty t^{-\frac{3}{2}} (f_j(x_0)-e^{t\Delta_D}f_j (x_0)) dt .
\]
We here notice that the integrand above is non-negative, 
since the $L^\infty$ norm of $e^{t\Delta_D}f_j $ is non-increasing. 
By a change of variable $t \mapsto 2^{-2j}t$, 
\[
\Lambda _D f_j(x_0) = c_1 2^j  \int_0^\infty t^{-\frac{3}{2}} (f_j(x_0)-e^{t2^{-2j}\Delta_D}f_j (x_0)) dt .
\]
We write 
\[
e^{t2^{-2j}\Delta_D}f_j = \Big( \sum_{l < j} + \sum _{l \geq j} \Big) \psi_l(\Lambda _D) e^{t2^{-2j}\Delta_D}f_j 
= \dfrac{\psi_j(\Lambda_D)}{1+2^{-2j} \Lambda_D ^2} e^{t2^{-2j}\Delta _D} f_j 
+ \sum _{l \geq j} \psi_l(\Lambda _D) e^{t2^{-2j}\Delta_D}f_j ,
\]
and apply the spectral multiplier theorem~\eqref{0901-1} to have that 
\[
\begin{split}
\left\| \dfrac{\psi_j(\Lambda_D)}{1+2^{-2j} \Lambda_D ^2} e^{t2^{-2j}\Delta _D } f_j  \right\|_{L^\infty} 
\leq C e^{-ct } \| f_j \|_{L^\infty} 
= C e^{-ct } f_j (x_0),
\end{split}
\]
and 
\[
\begin{split}
\left\| \sum _{l \geq j} \psi_l(\Lambda _D) e^{t2^{-2j}\Delta_D}f_j \right\|_{L^\infty} 
\leq 
& C \sum _{ l \geq j} \dfrac{2^{2j}}{2^{2l}} e^{-ct} \| f_j \|_{L^\infty} 
\leq C e^{-ct} f_j (x_0).
\end{split}
\]
We can then 
find a half line $[a,\infty)$ independent of $j$ such that 
\[
\| e^{t2^{-2j}\Delta_D}f_j \|_{L^\infty} 
\leq \frac{1}{2} f_j (x_0), \quad \text{ for } t \in [a,\infty). 
\]
We then obtain that 
\[
\Lambda _D f_j (x_0) \geq c_1 2^j  \int_a^\infty t^{-\frac{3}{2}} \left(f_j(x_0)-\frac{1}{2} f_j (x_0) \right) dt 
= \left( \frac{c_1}{2} \int_a^\infty t^{-\frac{3}{2}} dt \right) 2^j \| f_j \|_{L^\infty},
\]
which completes the proof. 
\end{pf}

\begin{cor}\label{cor:0902}
Let $\theta, f$ be smooth functions satisfying 
$\partial _t \theta + \Lambda_D \theta = f$. Then 
\begin{equation}\label{0906-11}
\partial _t \| \psi_j(\Lambda _D)\theta\|_{L^\infty} + c 2^j \| \psi_j(\Lambda _D) \theta \|_{L^\infty} 
\leq \| f \|_{L^\infty}. 
\end{equation}
\end{cor}
\begin{pf}
At a maximum point of $|\psi _j(\Lambda _D)f|$, 
we apply Lemma~3.2 in \cite{WZ-2011} (see also Lemma~2.2 in \cite{Iw-2020}) 
to the time derivative and Proposition~\ref{0906-3} to the fractional Laplacian, 
and obtain \eqref{0906-11}. 
\end{pf}

\subsection{Bilinear estimate}

\begin{prop}\label{prop:0830-10}
For every $f \in \dot B^1_{\infty,1} $ and $g \in \dot B^1_{\infty,1}$ 
\begin{equation}\label{0830-8}
\left\| (\nabla^{\perp}  f \cdot \nabla) g\right\|_{\dot B^0_{\infty,1}}
\leq C \| f \|_{B^1_{\infty,1}} \| g \|_{\dot B^1_{\infty,1}}. 
\end{equation}
Moreover, for every $f \in \dot B^0_{\infty,1}$ and $g \in \dot B^1_{\infty,1}$ 
\begin{equation}\label{0830-8-2}
\left\| (\nabla^{\perp} \Lambda _D ^{-1} f \cdot \nabla) g\right\|_{\dot B^0_{\infty,1}}
\leq C \| f \|_{\dot B^0_{\infty,1}} \| g \|_{\dot B^1_{\infty,1}}. 
\end{equation}
\end{prop}

\begin{pf}By the decomposition of the unity, we write
\[
f_j := \phi_j(\Lambda _D)f, g_j = \phi_j(\Lambda _D)g, 
\]
\[
\begin{split}
(\nabla^{\perp}  f \cdot \nabla) g
=
& \sum_{j,k,l \in \mathbb Z} \phi_j (\Lambda _D) 
\Big(  \nabla ^{\perp} f_k \cdot  \nabla g_l   \Big), 
\end{split}
\]
and we divide into the two cases when $j \geq \max \{ k,l\}$ and $j < \max\{ k,l \}$. 

We start by the case when $j \geq \max \{ k,l\}$. It follows from the term in the domain of the Dirichlet Lalacian 
(see~\eqref{0830-4}) that 
\[
\phi_j (\Lambda _D) 
\Big(  \nabla ^{\perp} f_k \cdot  \nabla g_l   \Big)
= \Lambda _D^{-2}\phi_j(\Lambda _D) 
 \left( (-\Delta _D) \Big(  \nabla ^{\perp} f_k \cdot  \nabla g_l   \Big)\right) .
\]
We then apply the Leibniz rule to $-\Delta _D$ and the inequalities in Lemmas~\ref{lem:0830-5}, \ref{lem:0830-6} to have 
\[
\begin{split}
\| \phi_j (\Lambda _D) 
\Big(  \nabla ^{\perp} f_k \cdot  \nabla g_l   \Big)
\|_{L^\infty} 
\leq & C 2^{-2j} (2^{3k+l} + 2^{k +3l} ) \| f_k \|_{L^\infty} \| g_l \|_{L^\infty}
\end{split}
\]
By the inequality above, we estimate 
\[
\begin{split}
& \sum _{j \geq \max\{ k,l \}}\| \phi_j (\Lambda _D) 
\Big(  \nabla ^{\perp} f_k \cdot  \nabla g_l   \Big)
\|_{L^\infty} 
\\
\leq & 
C  \sum _{j \geq \max\{ k,l \}} 2^{-2j} 2^{3k} \| f_k \|_{L^\infty} 2^l \| g_l \|_{L^\infty} 
+ C  \sum _{j \geq \max\{ k,l \}} 2^{-2j} 2^{k} \| f_k \|_{L^\infty} 2^{3l} \| g_l \|_{L^\infty} 
\\
\leq & 
C  \sum _{j\in \mathbb Z , k' \geq 0} 2^{-2j} 2^{3(j-k')} \| f_{j-k'} \|_{L^\infty} \| g \|_{\dot B^1_{\infty,1}}
+ C  \sum _{j \in \mathbb Z, l' \geq 0} 2^{-2j} \| f\|_{\dot B^1_{\infty,1}} 2^{3(j-l')} \| g_{j-l'} \|_{L^\infty} 
\\
\leq & 
C  \sum _{k' \geq 0} 2^{-2k'} \| f \|_{\dot B^1_{\infty,1}} \| g \|_{\dot B^1_{\infty,1}}
+ C  \sum _{ l' \geq 0}2^{-2l'}  \| f\|_{\dot B^1_{\infty,1}}\| g \|_{\dot B^1_{\infty,1}} .
\end{split}
\]

We next consider the case when $j < \max\{ k,l \}$, by dividing into two cases 
$k > l$ and $k \leq l$. When $k > l$, it follows from 
$ \nabla ^{\perp} f_k \cdot  \nabla g_l =\nabla ^\perp \big( f_k \nabla g_l \big) $ and \eqref{0830-2} that 
\[
\begin{split}
& 
\sum _{j < k , k > l}\| \phi_j (\Lambda _D) 
\Big(  \nabla ^{\perp} f_k \cdot  \nabla g_l   \Big)
\|_{L^\infty} 
\\
= & \sum _{j < k , k > l}\| \phi_j (\Lambda _D) 
\nabla ^{\perp} \cdot \big(   f_k \nabla g_l  \big)
\|_{L^\infty} 
\leq 
 C\sum _{j < k , k > l} 2^j \| f_k \|_{L^\infty} 2^l \| g_l \|_{L^\infty} 
\leq C\| f \|_{\dot B^1_{\infty,1} } \| g \|_{\dot B^1_{\infty,1}} .
\end{split}
\]
Analogously, 
\[
\begin{split}
& 
\sum _{j < l , k \leq l}\| \phi_j (\Lambda _D) 
\Big(  \nabla ^{\perp} f_k \cdot  \nabla g_l   \Big)
\|_{L^\infty} 
\\
= 
& \sum _{j < l , k \leq l}\| \phi_j (\Lambda _D) 
\nabla \Big( ( \nabla ^{\perp} f_k)   g_l   \Big)
\|_{L^\infty} 
\leq C \sum_{j < l , k \leq l} 2^j 2^k\| f_k \|_{L^\infty} \| g \|_{\dot B^1_{\infty,1}}
\leq C \|f\|_{\dot B^1_{\infty,1}} \| g \|_{\dot B^1_{\infty,1}} . 
\end{split}
\]
We obtain the first inequality \eqref{0830-8}. The second inequality follows from the 
lifting property~\eqref{0915-1}.   
\end{pf}

\subsection{Commutator estimates}

\begin{prop}\label{prop:0906-2}
Let $1 \leq q \leq \infty$ and  
\[
\psi_j (\lambda) := (1+2^{-2j-2}\lambda ^2)^{-1} - (1+2^{-2j}\lambda ^2)^{-1} 
= \dfrac{3}{4} \cdot \dfrac{2^{-2j} \lambda ^2}{(1+2^{-2j-2}\lambda^2) (1+2^{-2j}\lambda ^2)}. 
\]
Then there exists $C >0$ such that 
\begin{equation}\label{0902-11}
\begin{split}
\left\{\sum_{j \in \mathbb Z} 
\Big\| (\nabla ^\perp \Lambda _D^{-1} f \cdot \nabla) \psi_j(\Lambda _D) g 
- \psi_j(\Lambda_D) \Big( (\nabla ^\perp \Lambda _D^{-1} f \cdot \nabla) g \Big)
\Big\|_{L^\infty}^q
\right\}^{\frac{1}{q}}
\leq 
C \| f \|_{\dot B^{1}_{\infty,1}} \| g \|_{\dot B^{0}_{\infty,q}} ,
\end{split}
\end{equation}
\begin{equation}\label{0902-11-2}
\begin{split}
\left\{\sum_{j \in \mathbb Z} 
\Big\| (\nabla ^\perp \Lambda _D^{-1} f \cdot \nabla) \psi_j(\Lambda _D) g 
- \psi_j(\Lambda_D) \Big( (\nabla ^\perp \Lambda _D^{-1} f \cdot \nabla) g \Big)
\Big\|_{L^\infty}^q
\right\}^{\frac{1}{q}}
\leq 
C \| f \|_{\dot B^{\frac{1}{2}}_{\infty,\infty}} \| g \|_{\dot B^{\frac{1}{2}}_{\infty,q}} .
\end{split}
\end{equation}

\end{prop}
\begin{pf} We utilize the Littlewood-Paley dyadic decomposition $\{ \phi_j \}_{j \in \mathbb Z}$,  and the resolution of the identity 
\[
f = \sum_{k \in \mathbb Z} \phi_k(\Lambda _D)f =  \sum _{ k \in \mathbb Z} f_k , \quad 
g = \sum _{ l \in \mathbb Z} \phi_l (\Lambda_D )g = \sum _{l \in \mathbb Z} g_l. 
\]
We mainly discuss the proof of the second inequality \eqref{0902-11-2}, 
since the first one \eqref{0902-11} can be handled rather easier.

We write $f = S_j f + (1-S_j)f$ and start by an easier part having $(1-S_j) f$. 
The term with $(1-S_j)f$ does not require some cancellation due to 
commutator and we estimate two terms in the left hand side of \eqref{0902-11-2} separately. 
The first term in the left hand side of \eqref{0902-11-2} 
is estimated by using \eqref{0830-3} and \eqref{0830-2}
\[
\begin{split}
& 
\left\{ \sum_{j \in \mathbb Z} 
\Big\| (\nabla ^\perp \Lambda _D^{-1}(1-S_j) f \cdot \nabla) \psi_j(\Lambda _D) g \Big\|_{L^\infty}^q\right\}^{\frac{1}{q}}
\\
\leq 
&C \left\{ \sum_{j \in \mathbb Z} \Big( \sum_{k >j} 
  \| \nabla ^\perp \Lambda _D^{-1} f_k \|_{L^\infty} 
   \sum _{l \in \mathbb Z}\| \nabla \psi_j (\Lambda _D) g_l \|_{L^\infty}\Big) ^q
   \right\}^{\frac{1}{q}}
\\
\leq 
& C \left\{ \sum _{j \in \mathbb Z}  
 \Big( 2^{-\frac{1}{2}j} \Big( \sup_{k >j} 2^{\frac{1}{2}k} \| f_k \| _{L^\infty}\Big) 
 \sum_{l \in \mathbb Z} 2^l \frac{2^{-2j+ 2l}}{(1+2^{-2j+2l})^2}\| g_l \|_{L^\infty} 
 \Big) ^q \right\}^{\frac{1}{q}}
\\
\leq 
&  C \| f \|_{\dot B^\frac{1}{2}_{\infty,\infty}}  \left\{ \sum _{j \in \mathbb Z}  
 \Big( 
 \sum_{l \in \mathbb Z} 2^{\frac{1}{2}(l-j)} 2^{-2|j-l|}  \cdot 2^{\frac{1}{2}l} \| g_l \|_{L^\infty} 
 \Big) ^q \right\}^{\frac{1}{q}} 
\\
\leq 
& C \| f \|_{\dot B^\frac{1}{2}_{\infty,\infty}} \| g \|_{\dot B^{\frac{1}{2}}_{\infty,q}} .
\end{split}
\]
For the second term in the left hand side of \eqref{0902-11-2}, 
we decompose $g = S_j g + (1-S_j) g$ and have that for the term 
with $S_j g$
\[
\begin{split}
& 
\left\{ \sum_{j \in \mathbb Z} 
\Big\| \psi_j(\Lambda _D)  \Big( (\nabla ^\perp \Lambda _D^{-1}(1-S_j) f \cdot \nabla) 
  S_j g\Big) \Big\|_{L^\infty}^q\right\}^{\frac{1}{q}}
\\
\leq 
& C 
\left\{ \sum_{j \in \mathbb Z} 
\Big( \sum_{k > j} \| f_k \|_{L^\infty} \sum _{l \leq j} 2^l\| g_l \|_{L^\infty}
\Big) ^q\right\}^{\frac{1}{q}}
\\
\leq 
& C 
\left\{ \sum_{j \in \mathbb Z} 
\Big( 2^{-\frac{1}{2}j}\| f \|_{\dot B^{\frac{1}{2}}_{\infty,\infty}} \Big) ^q
  \Big( 2^{\frac{1}{4}j}\Big)^q \sum _{l \leq j} \Big(2^{\frac{3}{4}l}\| g_l \|_{L^\infty}
\Big) ^q\right\}^{\frac{1}{q}}
\\
\leq 
& C \| f \|_{\dot B^{\frac{1}{2}}_{\infty,\infty}} 
\left\{ \sum _{ l \in \mathbb Z}\sum_{j \geq l } 
\Big( 2^{-\frac{1}{4}j} \Big)^q  \Big(2^{\frac{3}{4}l}\| g_l \|_{L^\infty}
\Big) ^q\right\}^{\frac{1}{q}}
\leq C \| f \|_{\dot B^{\frac{1}{2}}_{\infty,\infty}} \| g \|_{\dot B^{\frac{1}{2}}_{\infty,q}} .
\end{split}
\]
For the term with $(1-S_j)g$, 
we write the divergence form, 
$(\nabla ^\perp \Lambda _D^{-1}(1-S_j) f \cdot \nabla) (1-S_j) g 
= \nabla \cdot \big( (\nabla ^\perp \Lambda _D^{-1}(1-S_j) f ) (1-S_j) g\big)$, 
and by \eqref{0830-2} 
\[
\begin{split}
& 
\left\{ \sum_{j \in \mathbb Z} 
\Big\| \psi_j(\Lambda _D)  
\nabla \cdot \big( (\nabla ^\perp \Lambda _D^{-1}(1-S_j) f ) (1-S_j) g\big)
 \Big\|_{L^\infty}^q\right\}^{\frac{1}{q}}
\\
\leq 
& 
\left\{ \sum_{j \in \mathbb Z} 
\Big\{ \sum _{j' \in \mathbb Z} \Big\| \psi_j (\Lambda _D)\phi_j(\Lambda _D)  
\nabla \cdot \big( (\nabla ^\perp \Lambda _D^{-1}(1-S_j) f ) (1-S_j) g\big)
 \Big\|_{L^\infty}\Big\} ^q\right\}^{\frac{1}{q}} 
\\
\leq 
& 
\left\{ \sum_{j \in \mathbb Z} 
\Big\{ \sum _{j' \in \mathbb Z} 2^{j'} 2^{-2|j-j'|}
\big\| (\nabla ^\perp \Lambda _D^{-1}(1-S_j) f ) (1-S_j) g
\big\|_{L^\infty}\Big\} ^q\right\}^{\frac{1}{q}}.
\end{split}
\]
We then have that 
\[
\begin{split}
& 
\left\{ \sum_{j \in \mathbb Z} 
\Big\{ \sum _{j' \in \mathbb Z} 2^{j'} 2^{-2|j-j'|}
\big\| (\nabla ^\perp \Lambda _D^{-1}(1-S_j) f ) (1-S_j) g
\big\|_{L^\infty}\Big\} ^q\right\}^{\frac{1}{q}}.
\\
\leq 
& C 
\left\{ \sum_{j \in \mathbb Z} 
\Big(\sum_{j' \in \mathbb Z} 2^{j'} 2^{-2|j-j'|}  \cdot 2^{-\frac{1}{2}j} 
 \| f \|_{\dot B^{\frac{1}{2}}_{\infty,\infty}} 
\sum _{l > j} \| g_l \|_{L^\infty}
\Big) ^q\right\}^{\frac{1}{q}}
\\
= 
& C \| f \|_{\dot B^{\frac{1}{2}}_{\infty,\infty}}  
\left\{ \sum_{j \in \mathbb Z} 
\Big( \sum _{j' \in \mathbb Z} 2^{(j'-j)} 2^{-2|j-j'|}
\sum _{l' > 0}  2^{-\frac{1}{2}l'} 2^{\frac{1}{2}(j+l')}\| g_{j+l'} \|_{L^\infty}
\Big) ^q\right\}^{\frac{1}{q}}
\\
\leq 
& C \| f \|_{\dot B^{\frac{1}{2}}_{\infty,\infty}} 
 \| g \|_{\dot B^{\frac{1}{2}}_{\infty,q}} .
\end{split}
\]
We next consdier the remainder part having $S_j f$, and in this case we need 
a cancellation due to the commutator. 
Since the first term in the left hand side belongs to the domain of the Dirichlet Lalacian, 
the Dirichlet Laplacian can act on it. 
 We then note that 
\[
\begin{split}
& (\nabla ^\perp \Lambda _D^{-1} S_j f \cdot \nabla) \frac{1}{1-2^{-2j}\Delta_D} g
\\
= & 
\frac{1}{1-2^{-2j}\Delta_D} (1-2^{-2j}\Delta) 
(\nabla ^\perp \Lambda _D^{-1}  S_j f \cdot \nabla) \frac{1}{1-2^{-2j}\Delta_D} g 
\\
=& 
\frac{1}{1-2^{-2j}\Delta_D} (\nabla ^\perp \Lambda _D^{-1}  S_j f \cdot \nabla) g
\\
& + \frac{-2^{-2j}}{1-2^{-2j}\Delta_D} 
\Big( (\nabla ^{\perp} \Lambda _D  S_j f \cdot \nabla) \frac{1}{1-2^{-2j}\Delta_D} g
     + (\nabla^\perp \nabla \Lambda _D^{-1}  S_j f \cdot \nabla\nabla ) 
        \frac{1}{1-2^{-2j}\Delta_D} g 
\Big) 
\\
=:& 
\frac{1}{1-2^{-2j}\Delta_D} (\nabla ^\perp \Lambda _D^{-1}  S_j f \cdot \nabla) g
+ R_j (f, g),
\end{split}
\]
and we can write 
\begin{gather}\notag 
 (\nabla ^\perp \Lambda _D^{-1} S_j f \cdot \nabla) \psi_j(\Lambda _D) g 
- \psi_j(\Lambda_D) (\nabla ^\perp \Lambda _D^{-1} S_jf \cdot \nabla) g  
= R_{j+1}(f, g) - R_j ( f, g ), 
\end{gather}
and it suffices to estimate 
\begin{equation}\notag 
\Big\{ \sum _{ j \in \mathbb Z} \|  R_j (f,g) \|_{L^\infty} ^q \Big\} ^{\frac{1}{q}}.
\end{equation}
We apply the boundedness of the resolvent \eqref{0902-12}, 
the boundedness of spectral multipliers with the first and the second derivatives 
\eqref{0830-2} and \eqref{0819-1} to have that 
\[
\begin{split}
& 
\| R_j (f,g) \|_{L^\infty} 
\\
\leq & 
C 2^{-2j} 
\left\{ \sum_{k \leq j} 2^{2k} \| f_k \|_{L^\infty} 
 \sum_{l \in \mathbb Z} \frac{2^l}{1+2^{-2j+2l}} \| g_l \|_{L^\infty}
+ \sum_{k \leq j} 2^k \| f_k \|_{L^\infty} 
   \sum_{l \in \mathbb Z} \frac{2^{2l}}{1+2^{-2j+2l}} \| g_l \|_{L^\infty}
 \right\}
\\
\leq & 
C \| f \|_{\dot B^{\frac{1}{2}}_{\infty,\infty}} 2^{-2j}  
\left\{ 2^{\frac{3}{2}j} \cdot 
 \sum_{l \in \mathbb Z} 2^{\frac{1}{2}j} \cdot \frac{2^{\frac{1}{2}(l-j)}}{1+2^{-2j+2l}}2^{\frac{1}{2}l}\| g_l \|_{L^\infty}
+ 2^{\frac{1}{2}j}
   \sum_{l \in \mathbb Z} 2^{\frac{3}{2}j} \frac{2^{\frac{3}{2}(l-j)}}{1+2^{-2j+2l}} 2^{\frac{1}{2}l}\| g_l \|_{L^\infty}
 \right\}
\\
\leq & 
C \| f \|_{\dot B^{\frac{1}{2}}_{\infty,\infty}} 
\left\{ 
 \sum_{l \in \mathbb Z} 2^{-\frac{1}{2}|l-j|}\| g_l \|_{L^\infty}
+ 
   \sum_{l \in \mathbb Z} 2^{-\frac{3}{2}|l-j|} \| g_l \|_{L^\infty}
 \right\}
\leq 
C \| f \|_{\dot B^{\frac{1}{2}}_{\infty,\infty}} 
 \sum_{l \in \mathbb Z} 2^{-\frac{1}{2}|l-j|}\| g_l \|_{L^\infty}, 
\end{split}
\]
By taking the $\ell^q$ norm of the term above, we obtain \eqref{0902-11-2}.  
We argue for the first inequality \eqref{0902-11} here, since this part is crucial. 
Similarly to the above, 
\[
\begin{split}
& 
\| R_j (f,g) \|_{L^\infty} 
\\
\leq & 
C \| f \|_{\dot B^{1}_{\infty,1}} 2^{-2j}  
\left\{ 2^{j} \cdot 
 \sum_{l \in \mathbb Z}  \frac{2^{l}}{1+2^{-2j+2l}}\| g_l \|_{L^\infty}
+
   \sum_{l \in \mathbb Z}  \frac{2^{2l }}{1+2^{-2j+2l}} \| g_l \|_{L^\infty}
 \right\}
\\
\leq 
& C \| f \|_{\dot B^{\frac{1}{2}}_{\infty,\infty}} 
 \sum_{l \in \mathbb Z} \Big( 2^{-|l-j|} + 2^{-2|l-j|} \Big) \| g_l \|_{L^\infty}.
\end{split}
\]
and by taking the $\ell^1$ norm and the Young inequality,  we conclude \eqref{0902-11}. 
\end{pf}

\section{Proof of Theorem \ref{thm:1}}

Let the initial data $\theta_0 \in \dot B^0_{\infty ,1}$. Let $\{ \theta _n \}_{n=1}^\infty$ be defined by for $n = 1$
\begin{equation}\label{0902-4-0}
\begin{cases}
\partial _{t} \theta_{1} + \Lambda _D \theta_{1}  = 0 ,
\\
\theta _{1}(0) = \displaystyle \sum _{ j \leq 1 } \phi_j (\Lambda _D) \theta _0 ,
\end{cases}
\end{equation}
and for $n = 2,3, \cdots$ 
\begin{equation}\label{0902-4}
\begin{cases}
\partial _{t} \theta_{n} + \Lambda _D \theta_{n} 
 \displaystyle 
 + (\nabla ^\perp \Lambda _D^{-1} \theta _{n-1} \cdot \nabla) \theta _{n} 
 = 0 ,
\\
\theta _{n}(0) = \displaystyle \sum _{ j \leq n } \phi_j (\Lambda _D) \theta _0 =: \theta_{0,n} ,
\end{cases}
\end{equation}
It is easy to see that $\theta _1$ is well-defined, since it is a linear solution and we need to 
prove the existence of $\theta _n$ for $n \geq 2$ for given $\theta _{n-1}$. 

\begin{prop}\label{prop:0906-1}
For every $\theta _0 \in \dot B^0_{\infty,1}$, 
$T > 0$ exists such that there exists a function 
$\theta _n \in C([0,T) , \dot B^0_{\infty,1} ) \cap L^1 (0,T ; \dot B^1_{\infty,1})$ 
satisfying 
$\partial _t \theta _n \in L^1 (0,T ; \dot B^1_{\infty,1})$ 
and the equation \eqref{0902-4-0} for $n = 1$, 
the equation \eqref{0902-4} for $n \geq 2$. 
\end{prop}

\begin{pf}
The case when $n = 1$ follows obviously due to maximum regularity estimate \eqref{0902-2} in the time interval $[0,\infty)$, and let $n \geq 2$ and 
we assume that a solution $\theta _{n-1} $ exists such that 
$\theta _{n-1} \in C([0,T) , \dot B^0_{\infty,1} ) \cap L^1 (0,T ; \dot B^1_{\infty,1})$ 
satisfying 
$\partial _t \theta _{n-1} \in L^1 (0,T ; \dot B^1_{\infty,1})$, 
where the existence time $T$ will be disscussed later to be 
independent of $n$. 

We now show the existence of $\theta_n$. It is possible to obtain a local solution, where 
the existence time depends on $n$ by Galerkin approximations (see~\cite{CoIg-2017,CoNg-2018-2}), 
but we give a self-contained proof to make the independency for the existence time clear in our framework. 
To this end, we approximate 
solutions by $\theta _\varepsilon$, a solution of the following equation.
\begin{equation}\label{0902-5}
\begin{cases}
\partial _{t} \theta_{\varepsilon} + \Lambda _D \theta_{\varepsilon} -\varepsilon \Delta \theta _\varepsilon
 + (\nabla ^\perp \Lambda _D^{-1} \theta _{n-1} \cdot \nabla) \theta _{\varepsilon} 
 = 0 ,
\\
\theta _{n}(0) =  \theta_{0,n} ,
\end{cases}
\end{equation}
where $\varepsilon > 0$. 
We construct a solution $\theta _\varepsilon \in C([0,T], H^2) 
\cap C^1 ([0,T] , H^1)$ of \eqref{0902-5}, 
to obtain a solution $\theta_n$ of \eqref{0902-4} 
by passing to the limit as $\varepsilon \to 0$. 
We mainly discuss in the case when $n = 2$ with several estimates possible to be applied 
to the cases when $n \geq 3$.

Step 1 (Solution in a short interval $[0,T_\varepsilon]$ when $n = 2$). 
We consider the integral equation 
\begin{equation}\label{0902-9}
\theta_\varepsilon (t) 
= e^{-t(\Lambda _D+\varepsilon \Lambda _D^2)} \theta_{0,n} 
+ \int_0^t e^{-(t-\tau)(\Lambda _D+\varepsilon \Lambda _D^2)}
\Big( (\nabla ^\perp \Lambda _D^{-1} \theta _{n-1} \cdot \nabla) \theta _{\varepsilon} 
\Big) d\tau , 
\quad \text{with }n= 2 .
\end{equation}
It is easy to check that 
\[
\| \theta_{0,n} \|_{H^2} \leq C 2^{2n} \| \theta_0 \|_{L^2} 
\leq C 2^{2n} \| \theta_0 \|_{L^\infty} < \infty .
\]
Also, it is not difficult to show the following inequality
\[
\begin{split}
& 
\left\| 
\int_0^t e^{-(t-\tau)(\Lambda _D+\varepsilon \Lambda _D^2)}
\Big( (\nabla ^\perp \Lambda _D^{-1} f \cdot \nabla) g
\Big) d\tau
\right\|_{H^2}
\leq 
\dfrac{CT^{\frac{1}{2}}}{\varepsilon^{\frac{1}{2}}} \| f \|_{L^\infty (0,T ; H^2)} \| g \|_{L^\infty (0,T ; H^2)} .
\end{split}
\]
In fact, 
\[
\| e^{-(t-\tau)(\Lambda _D+\varepsilon \lambda _D^2)} \|_{H^1 \to H^2} 
\leq C \big( (t-\tau) \varepsilon \big)^{-\frac{1}{2}} , 
\]
and by the H\"older inequality, \eqref{0902-6}, \eqref{0902-8}, \eqref{0902-7}, \eqref{0819-1-2}
\[
\begin{split}
\| (\nabla ^\perp \Lambda _D^{-1} f \cdot \nabla) g \|_{H^1}
\leq  
& C \big( \| \nabla \nabla ^\perp \Lambda _D^{-1} f \|_{L^4} \| \nabla g \|_{L^4}
+  \| \nabla ^\perp \Lambda _D^{-1}f \|_{L^\infty} \| \nabla ^2 g \|_{L^2} \big)
\\
\leq & C \big( \| f \|_{\dot B^1_{4,1}} \| g \|_{\dot B^1_{4,1}} 
 + \| f \|_{\dot B^0_{\infty,1}} \| g \|_{H^2} \big) 
\\
\leq& C \| f \|_{\dot B^2_{2,2}} \| g \|_{\dot B^2_{2,2}} ,
\end{split}
\]
which prove the inequality above. We can then apply the Banach fixed point theorem 
to obtain a solution $\theta _\varepsilon \in C([0,T _\varepsilon], H^2 )$ 
of \eqref{0902-9}, 
and it satisfies $\partial_t \theta _{\varepsilon} \in C^1 ([0,T_\varepsilon], H^1)$, 
where $T_\varepsilon \leq \varepsilon/ (C\| \theta _{0,n} \|_{H^2})^2$, 
with $n = 2$. 

Step 2 (A priori estimate when $n = 2$). 
Suppose that $\theta _\varepsilon$ is a solution of \eqref{0902-9} in the time interval $[0, \infty)$.  
For sufficiently small $\delta > 0$, we prove that  there exists $T > 0$ independent of $\varepsilon$ such that 
\[
\| \theta_\varepsilon \|_{L^1 (0,T ; \dot B^1_{\infty,1})} \leq 2 \delta .
\]
$\psi_j(\Lambda_D)$ acting on the equation, we write 
\begin{equation}\label{0906-1}
\begin{split}
& \partial _t \psi_j(\Lambda _D) \theta _\varepsilon 
+ \Lambda _D \psi_j(\Lambda _D) \theta _\varepsilon  
- \varepsilon \Delta _D \psi_j(\Lambda _D) \theta _\varepsilon  
+ (\nabla ^\perp \Lambda _D^{-1} \theta_{n-1} \cdot \nabla) \psi_j(\Lambda _D) \theta _{\varepsilon} 
\\
=
&  (\nabla ^\perp \Lambda _D^{-1} \theta_{n-1} \cdot \nabla) \psi_j(\Lambda _D) \theta _{\varepsilon} 
-\psi_j(\Lambda _D) 
 \Big( (\nabla ^\perp \Lambda _D^{-1} \theta_{n-1} \cdot \nabla) \theta _{\varepsilon} 
 \Big)
\\
=& 
\Big[\nabla ^\perp \Lambda _D^{-1} \theta_{n-1} \cdot \nabla , \,\, \psi_j(\Lambda _D) 
  \Big] \theta_\varepsilon .
\end{split}
\end{equation}
For almost every $t$, at a maximum point of $|\psi_j (\Lambda _D)\theta_\varepsilon|$, 
it follows from $-\Delta _D \psi_j (\Lambda _D)\theta_\varepsilon $ having 
the same sign as $\psi_j (\Lambda _D)\theta_\varepsilon$ and Corollary~\ref{cor:0902} that 
\[
\partial _t \| \psi_j(\lambda_D)\theta_{\varepsilon} \|_{L^\infty } 
+ c 2^j \| \psi_j(\Lambda _D) \theta _\varepsilon\|_{L^\infty} 
\leq \Big\| \Big[\nabla ^\perp \Lambda _D^{-1} \theta_{n-1} \cdot \nabla , \,\, \psi_j(\Lambda _D) 
  \Big] \theta_\varepsilon 
     \Big \|_{L^\infty}  ,
\]
which yields that 
\begin{equation}\label{0902-10}
\| \psi_j(\Lambda _D) \theta_\varepsilon \|_{L^\infty} 
\leq e^{-ct 2^j} \| \psi_j(\Lambda_D) \theta_0 \|_{L^\infty} 
+ \int_0^t e^{-c(t-\tau)2^j} 
\Big\| \Big[\nabla ^\perp \Lambda _D^{-1} \theta_{n-1} \cdot \nabla , \,\, \psi_j(\Lambda _D) 
  \Big] \theta_\varepsilon 
     \Big \|_{L^\infty} 
 ~d\tau  .
\end{equation}
By summing over $j \in \mathbb Z$ and applying \eqref{0902-11-2}, we have 
\begin{equation}\label{0903-8}
\begin{split}
\| \theta(t) \|_{\dot B^0_{\infty,1}} 
\leq 
& C\| \theta_0 \|_{\dot B^0_{\infty,1}} 
+ C\sum_{j \in \mathbb Z } \int_0^t \Big\| \Big[\nabla ^\perp \Lambda _D^{-1} \theta_{n-1} \cdot \nabla , \,\, \psi_j(\Lambda _D) 
  \Big] \theta_\varepsilon 
     \Big \|_{L^\infty} 
  d\tau 
\\
\leq & 
C \| \theta_0 \|_{\dot B^0_{\infty,1}} 
+ C \int_0^t \| \theta_{n-1} \|_{\dot B^\frac{1}{2}_{\infty,1}} \| \theta_\varepsilon \|_{\dot B^\frac{1}{2}_{\infty,1}} d\tau 
\\
\leq & 
C \| \theta_0 \|_{\dot B^0_{\infty,1}} 
+ C \| \theta_{n-1} \|_{L^2 (0,t;\dot B^\frac{1}{2}_{\infty,1})} 
    \| \theta_\varepsilon \|_{L^2 (0,t; \dot B^\frac{1}{2}_{\infty,1})} .
\end{split}
\end{equation}
Also, by multiplying \eqref{0902-10} by $2^j$ and taking the $L^1 $ norm for time variable and 
the $\ell ^1$ norm 
\begin{equation}\label{0903-9}
\begin{split}
\| \theta_\varepsilon \|_{L^1(0,T ; \dot B^1_{\infty,1})} 
\leq 
& C \| e^{-ct\Lambda _D}\theta_0 \|_{L^1 (0,T ; \dot B^1_{\infty,1})}
+ C \int_0^T \| \theta_{n-1} \|_{\dot B^\frac{1}{2}_{\infty,1}} \| \theta_\varepsilon \|_{\dot B^\frac{1}{2}_{\infty,1}} d\tau 
\\
\leq & 
C \| e^{-ct\Lambda _D}\theta_0 \|_{L^1 (0,T ; \dot B^1_{\infty,1})}
+ C \| \theta_{n-1} \|_{L^2 (0,T;\dot B^\frac{1}{2}_{\infty,1})} 
    \| \theta_\varepsilon \|_{L^2 (0,T; \dot B^\frac{1}{2}_{\infty,1})} .
\end{split}
\end{equation}
We then take $\delta > 0$ such that 
\[
\delta + C (2\delta)^2 \leq 2 \delta ,
\]
and $T > 0$ such that 
\begin{equation}\label{0903-10}
C \max \left\{ \| e^{-ct\Lambda _D}\theta_0 \|_{L^1 (0,T ; \dot B^1_{\infty,1})}, 
 \Big( \| \theta_0 \|_{\dot B^0_{\infty,1}}+  (2\delta)^2 \Big) ^{\frac{1}{2}}  
  \| e^{-t\Lambda _D } \theta _0  \|_{L^1 (0,T ; \dot B^{1}_{\infty,1})} ^{\frac{1}{2}}, 
   ( 2 \delta )^2 
 \right\} \leq \delta .  
\end{equation}
We have from \eqref{0903-9} that 
\begin{equation}\label{0903-11}
\| \theta_\varepsilon \|_{L^1(0,T ; \dot B^1_{\infty,1})} \leq \delta  + C (2\delta )^2 
\leq 2 \delta  , \quad \text{ with } n = 2.
\end{equation}
Step 3 (Independent existence time of $\varepsilon $ when $n = 2$). 
By Step 1 above, it is sufficient to have a boundedness in $H^2$ 
independent of $\varepsilon$, $n$. 
$\Lambda _D^2 = -\Delta_D $ acting on the equation \eqref{0902-5}, we have 
\[
\partial_t \Lambda _D^2 \theta _\varepsilon + \Lambda (\Lambda _D^2 \theta_\varepsilon)
+ (\nabla ^{\perp} \Lambda _D^{-1} \theta_{n-1} \cdot \nabla) \Lambda _D^2 \theta_\varepsilon
= 2(\nabla \nabla ^{\perp} \Lambda _D^{-1} \theta_{n-1} \cdot \nabla)  \nabla \theta_\varepsilon
+ (\nabla ^{\perp} \Lambda _D \theta_{n-1} \cdot \nabla)  \theta_\varepsilon.
\]
Multiplication by $\Lambda _D^2 \theta_\varepsilon$, integration over the domain 
and the H\"older inequality give 
\[
\begin{split}
\frac{1}{2} \partial _t \| \theta_\varepsilon \|_{H^2}^2 
\leq & 
\frac{1}{2} \partial _t \| \theta_\varepsilon \|_{H^2}^2 + \| \theta_{\varepsilon} \|_{H^{\frac{5}{2}}}^2
\\
\leq 
&\| \nabla \nabla ^{\perp} \Lambda _D^{-1} \theta_{n-1} \|_{L^\infty} 
  \| \nabla ^2 \theta_\varepsilon \|_{L^2} \| \Lambda _D ^2 \theta_{\varepsilon}\|_{L^2}
 + \| \nabla ^\perp \Lambda _D \theta_{n-1} \|_{L^2} \| \nabla \theta_\varepsilon \|_{L^\infty} 
   \| \Lambda _D^2 \theta_\varepsilon \|_{L^2}
\\
\leq 
& C \| \theta_{n-1} \|_{\dot B^1_{\infty,1}} \| \theta_\varepsilon \|_{H^2} ^2 
+ C \| \theta_{n-1} \|_{H^2} \| \theta_\varepsilon \|_{\dot B^1_{\infty,1}} \| \theta_\varepsilon \|_{H^2}. 
\end{split}
\]
By the Young inequality and integrating over a time interval, when $0 \leq t \leq T$, 
\[
\begin{split}
& 
\| \theta_\varepsilon \|_{H^2}^2
\\
\leq 
& \| \theta_0 \|_{H^2}^2 
+ C \int_0^t \Big(  \| \theta_{n-1} \|_{\dot B^1_{\infty,1}} \| \theta_\varepsilon \|_{H^2} ^2 
     +  \| \theta_\varepsilon \|_{\dot B^1_{\infty,1}} \| \theta_{n-1} \|_{H^2}^2 
     +  \| \theta_\varepsilon \|_{\dot B^1_{\infty,1}} \| \theta_\varepsilon \|_{H^2}^2
     \Big) d\tau 
\\
\leq 
& \| \theta_0 \|_{H^2}^2 
 + C \| \theta_{\varepsilon} \|_{L^1 (0,T;\dot B^1)} \| \theta_{n-1} \|_{L^\infty(0,T ; H^2)}^2
 + C \int _0^t 
  \Big(  \| \theta_{n-1} \|_{\dot B^1_{\infty,1}} + \| \theta_{\varepsilon} \|_{\dot B^1_{\infty,1}}
  \Big)  \| \theta_\varepsilon \|_{H^2} ^2 d\tau ,
\end{split}
\]
and the Gronwall inequality implies that 
\[
\begin{split}
& \| \theta_\varepsilon (t) \|_{H^2}^2 
\\
\leq 
& \Big( \| \theta_0 \|_{H^2}^2 
       + C \| \theta_{\varepsilon} \|_{L^1 (0,T;\dot B^1_{\infty ,1})} \| \theta_{n-1} \|_{L^\infty(0,T ; H^2)}^2
\Big) 
\exp\left\{C \int_0^t  
    \Big(  \| \theta_{n-1} \|_{\dot B^1_{\infty,1}} + \| \theta_{\varepsilon} \|_{\dot B^1_{\infty,1}}
  \Big) d\tau \right\} .
\end{split}
\]
For $T > 0$ defined by \eqref{0903-10}, we have from \eqref{0903-11} that 
\[
\| \theta_\varepsilon (t) \|_{H^2}^2 
\leq 
\Big( \| \theta_0 \|_{H^2}^2 
       + C \cdot 2 \delta \cdot (2 \delta)^2 
\Big) 
\exp\left\{C (2 \delta + 2 \delta ) \right\}, \quad 
\text{ as long as } 0 < t \leq T,
\]
which yields the existence of $\theta _\varepsilon \in C([0,T] , H^2)$ for all $\varepsilon > 0$.

Step 4 (Convergence as $\varepsilon \to 0$ and the existence of $\theta_n$ when $n = 2$).  
Let $0 < \varepsilon' < \varepsilon$ and the difference of 
$\theta _\varepsilon - \theta_{ \varepsilon'}$ satisfies 
\[
\partial _t (\theta_\varepsilon - \theta_{\varepsilon'}) 
+ \Lambda _D (\theta_\varepsilon - \theta_{\varepsilon'}) 
- \varepsilon \Delta _D ( \theta _\varepsilon - \theta_{\varepsilon'} )
- (\varepsilon - \varepsilon ') \Delta _D \theta _\varepsilon' 
+ (\nabla ^\perp \Lambda _D ^{-1} \theta _{n-1} \cdot \nabla) (\theta_\varepsilon - \theta_{\varepsilon'})
= 0 .
\]
The $L^2$ inner product with $-\Delta _D(\theta_\varepsilon - \theta_{\varepsilon'})$ gives 
\[
\begin{split}
\frac{1}{2}\| \theta_\varepsilon(t) - \theta_{\varepsilon'} (t)\|_{H^1}^2 
\leq 
& \int_0^t 
\Big( (\varepsilon - \varepsilon') \| \theta_\varepsilon' \|_{H^2} 
  \| \theta_\varepsilon - \theta_{\varepsilon'} \|_{H^2}
  + \| \nabla \nabla ^\perp \Lambda _D^{-1}\theta_{n-1} \|_{L^\infty} 
    \| \theta_\varepsilon - \theta_{\varepsilon'} \|_{H^1} ^2
    \Big)
  d\tau 
\\
\leq 
& 2 (\varepsilon - \varepsilon ') T  \sup_{t \in [0,T] } 
  \| \theta_\varepsilon \|_{H^2} \| \theta_{\varepsilon '} \|_{H^2}
  +C  \int _0^t \| \theta_{n-1} \|_{\dot B^1_{\infty,1}} \| \theta_\varepsilon - \theta_{\varepsilon'} \|_{H^2}^2 d\tau .
\end{split}
\]
The Gronwall inequality implies that 
\[
\begin{split}
\sup _{t \in [0,T]}\| \theta_\varepsilon(t) - \theta_{\varepsilon'} (t)\|_{H^1}^2 
\leq 
&4 (\varepsilon - \varepsilon ') T  \sup_{t \in [0,T] } 
  \exp\left\{ C \| \theta_{n-1} \|_{L^1(0,T ; \dot B^1_{\infty,1})}\right\}
\\
\to 
& 0 \quad \text{ as } \varepsilon, \varepsilon ' \to 0 ,
\end{split}
\]
which implies that $\{ \theta_\varepsilon\}_{\varepsilon }$ satisfies a Cauchy condition, 
and 
we then obtain a limit function in 
$C([0,T], H^1 )$ and the uniform boundedness yields it also belongs to 
$L^\infty (0,T ; H^2 ) \cap L^\infty (0,T;\dot B^0_{\infty, 1}) 
\cap L^1 (0,T ; \dot B^1_{\infty,1})$ and the time derivative is in 
$L^1 (0,T \dot B^0_{\infty,1})$.  
By taking $\varepsilon \to 0$ for the integral equation \eqref{0902-9} 
in the topology of $L^\infty (0,T; L^2)$, 
we have that the limit function is a solution of \eqref{0902-4}, 
and hence we obtain $\theta_n$ with $n = 2$.

Step 5 (The case when $n \geq 3$). 
We use an induction argument. It is possible to argue analogously to Steps 1, 2, 3, 4 
and we obtain $\theta_n$ 
under a assumption that $\theta _{n-1}$ exists as 
a solution of \eqref{0902-4} with $n$ replaced by $n-1$ such that 
it belongs to $C([0,T], H^2) \cap L^\infty (0,T ; \dot B^0_{\infty,1})\cap 
L^1 (0,T, \dot B^1 _{\infty,1})$, where the existence time $T$ is same as 
Step 3. 
\end{pf}

\begin{prop}\label{prop:0906-2}
Let $\theta_n$ be obtained in Proposition~\ref{prop:0906-1}. 
Then there exists $T_0 \leq T$ such that 
$\theta_n$ converges in 
$L^\infty (0,T_0 ; \dot B^0_{\infty,1}) \cap L^1 (0,T_0;\dot B^1_{\infty,1})$ 
as $n \to \infty$, 
and the limt function, $\theta$, is a unique solution of \eqref{QG1} with the initial 
condition \eqref{QG2} in 
$L^\infty (0,T_0 ; \dot B^0_{\infty,1}) \cap L^1 (0,T_0;\dot B^1_{\infty,1})$. 
\end{prop}

\begin{pf}
Similarly to \eqref{0906-1}, we write 
\begin{equation}\label{0906-3}
\begin{split}
& \partial _t \psi_j(\Lambda _D) (\theta _{n+1} - \theta_n )
+ \Lambda _D \psi_j(\Lambda _D) (\theta _{n+1} - \theta_n )
+ (\nabla ^\perp \Lambda _D^{-1} \theta_{n} \cdot \nabla) \psi_j(\Lambda _D) 
   (\theta _{n+1} - \theta_n )
\\
&= 
\Big[\nabla ^\perp \Lambda _D^{-1} \theta_{n} \cdot \nabla , \,\, \psi_j(\Lambda _D) 
  \Big] (\theta _{n+1} - \theta_n ) 
- \psi_j (\Lambda _D) 
  \Big( \big( \nabla ^{\perp} \Lambda _D^{-1} (\theta_n - \theta_{n-1}) \cdot \nabla 
  \big) \theta_n 
  \Big) .
\end{split}
\end{equation}
We argue similarly to \eqref{0902-10}, \eqref{0903-8}, \eqref{0903-9}, 
but 
we apply a commutator estimate \eqref{0902-11} for the first term in the right 
hand side of the equality above instead of \eqref{0902-11-2}, 
and a bilinear estimate \eqref{0830-8-2} for the second term. 
\[
\begin{split}
& \| \theta_{n+1} - \theta_n 
\|_{L^\infty(0,t ; \dot B^0_{\infty,1}) \cap L^1 (0,t ; \dot B^1_{\infty,1})} 
\\
\leq
& C \| \phi_{n+1} (\Lambda _D)\theta_0 \|_{L^\infty}  
+C \int_0^t \| \theta_{n} \|_{\dot B^1_{\infty,1}}
   \| \theta_{n+1} - \theta_n \|_{\dot B^0_{\infty,1}} d\tau 
+ C \int_0^t \| \theta_n - \theta_{n-1} \|_{\dot B^0_{\infty,1}}  
      \| \theta_n \|_{\dot B^1_{\infty,1}} d\tau .
\end{split}
\]
We here write 
\[
D_{n+1}(t) := \| \theta_{n+1} - \theta_n 
\|_{L^\infty(0,t ; \dot B^0_{\infty,1}) \cap L^1 (0,t ; \dot B^1_{\infty,1})} 
\]
and then obtain by the inequality above that 
\[
D_n(t) 
\leq C \| \phi_{n+1} (\Lambda _D)\theta_0 \|_{L^\infty}  
+ C \| \theta_n \|_{L^1 (0,t ; \dot B^1_{\infty,1})} 
  (D_{n+1}(t) + D_n (t)). 
\]
We may have 
\[
C \| \theta_n \|_{L^1 (0,T ; \dot B^1_{\infty,1})} \leq \frac{1}{2}
\]
by taking the time interval shorter if necessary, and it yields that 
\[
D_{n+1}(t) \leq 2 C \| \phi_{n+1} (\Lambda _D)\theta_0 \|_{L^\infty}  
+ \frac{1}{2} D_n(t).
\]
We deduce that 
\[
\lim_{n\to\infty}\sum_{k=2}^n (\theta_k - \theta_{k-1} ) + \theta _1
\quad \text{ in } 
L^\infty (0,T ; \dot B^0_{\infty,1}) \cap L^1 (0,T ; \dot B^1_{\infty,1})
\]
exists, and let us define $\theta$ by the limit above. 
We can also prove that $\theta$ belongs to 
$C([0,T ]; \dot B^0_{\infty,1})$ and satisfies 
$\partial _t \theta \in L^1 (0,T ; \dot B^0_{\infty,1})$ 
and the equation $\partial _t \theta + \Lambda _D \theta + 
(\nabla ^\perp \Lambda_D^{-1} \theta \cdot \nabla) \theta = 0$. 

The uniqueness of solutions follows from an argument starting by 
a similar equality to \eqref{0906-3} for the difference of two solutions $\theta, \widetilde \theta $
with the same data, which implies 
\[
\begin{split}
& 
\| \theta -\widetilde \theta\|_{L^\infty (0,T ; \dot B^0_{\infty,1}) \cap L^1 (0,T ; \dot B^1_{\infty,1})}
\\
\leq 
& 
C \big( \| \theta \|_{L^1(0,T; \dot B^1_{\infty,1})} + \| \widetilde \theta \|_{L^1(0,T; \dot B^1_{\infty,1})}
\big)
\| \theta -\widetilde \theta\|_{L^\infty (0,T ; \dot B^0_{\infty,1}) \cap L^1 (0,T ; \dot B^1_{\infty,1})} .
\end{split}
\]
Therefore, we conclude the existence and the uniqueness of the local solution. 
\end{pf}

\bigskip 

\noindent 
{\bf Proof of Theorem~\ref{thm:1}}. 
Local existence of solutions follows from Proposition~\ref{prop:0906-2}, 
and Lemma~\ref{lem:0906-4} proves the zero boundary value. 
We also see that global existence for small data holds, since the constants,  
$C$, appearing in the proofs of Propositions~\ref{prop:0906-1}, \ref{prop:0906-2} 
are independent of the time interval and we can have the global existence result. 
\hfill $\Box$

\bigskip 

\noindent 
{\bf Remark}. 
We here introduce an alternative simple proof for the intepolation index $q = 1$ 
by using the integral equation, 
seeking a fixed point of the integral equation. 
\[
\theta (t) = e^{-t\Lambda _D} \theta_0 
-\int_0^t e^{-(t-\tau )\Lambda _D} 
  \Big( \big(\nabla ^\perp \Lambda _D ^{-1} \theta \cdot \nabla \big)
    \theta \Big) d\tau .
\]
Let 
\[
\Psi(\theta):= e^{-t\Lambda _D} \theta_{0,n} - \int_0^t e^{-(t-\tau) \Lambda _D} \Big( \big(\nabla ^\perp \Lambda _D ^{-1} \theta \cdot \nabla \big)
    \theta \Big) d\tau ,
\]
\[
X := L^\infty(0,\infty ; \dot B^0_{\infty,1} ) \cap L^1 (0,\infty;\dot B^1_{\infty,1}).
\]
We then have from maximum regularity estimate \eqref{0902-2}, \eqref{0902-3} 
and the bilinear estimate \eqref{0830-8-2} that  
\[
\begin{split}
 \| \Psi(\theta) \|_{X} 
\leq 
& 
C \| \theta_0 \|_{\dot B^0_{\infty,1}} 
+ C \int_0^\infty \Big\| \Big( \big(\nabla ^\perp \Lambda _D ^{-1} \theta \cdot \nabla \big)
    \theta \Big) \Big\|_{\dot B^0_{\infty,1}} d\tau 
\\
\leq 
& C \| \theta_0 \|_{\dot B^0_{\infty,1}} 
+ C \| \theta \|_{X}^2  ,
\\
 \| \Psi(\theta) - \Psi(\widetilde \theta) \|_{X} 
 \leq 
& C ( \| \theta \|_{X} + \| \widetilde \theta \|  _{X } )
    \| \theta - \widetilde \theta \|_{X},
\end{split}
\]
which allows to apply contraction argument. We conclude the global existence for small 
initial data.

\appendix 

\section{Equivalency of two resolutions} \label{Appen_1}

\noindent 
{\bf Proof of Lemma~\ref{lem:0906-8}}. 
The resolution with $\{ \phi_j \}_{j \in \mathbb Z}$ follows from 
Lemma~4.5 in \cite{IMT-2019}. 
We write 
\[
\| f \|_{\widetilde {\dot B^s_{p,q}}} 
:= \left\{ \sum _{ j \in \mathbb Z} 
   \Big( 2^{sj} \| \psi_j(\Lambda _D) f \|_{L^p } \Big) ^q
   \right\}^{\frac{1}{q}} ,
\]
and show the equivalency of the norm defined by $\{ \psi _j \}_{j \in \mathbb Z}$, 
where we see
\[
\psi_j (\lambda) = 
\dfrac{3}{4} \cdot \dfrac{2^{-2j} \lambda ^2}{(1+2^{-2j-2}\lambda^2) (1+2^{-2j}\lambda ^2)}. 
\]
Let $\Phi_j := \phi_{j-1} + \phi_j + \phi_{j+1}$, which satisfies 
$\phi_j = \Phi_j \phi_j$. 
It follows from the resolution by $\{ \phi_{j'} \}_{j' \in \mathbb Z}$ 
and the boundedness of the spectral multiplier that 
\[
\begin{split}
\| f \|_{\widetilde{\dot B^s_{p,q}}} 
\leq 
& \left\{ \sum _{ j \in \mathbb Z} 
   \Big( \sum _{j' \in \mathbb Z} 2^{sj} 
     \| \psi_j(\Lambda _D) \Phi_{j'}(\Lambda _D) \phi_{j'} (\Lambda _D)f \|_{L^p } \Big) ^q
   \right\}^{\frac{1}{q}} 
\\
\leq 
& C \left\{ \sum _{ j \in \mathbb Z} 
   \Big( \sum _{j' \in \mathbb Z} 2^{sj} 
    \frac{2^{-2j+2j'}}{(1+2^{-2j + 2j')^2}}
     \| \phi_{j'} (\Lambda _D)f \|_{L^p } \Big) ^q
   \right\}^{\frac{1}{q}} 
\\
\leq 
& C \left\{ \sum _{ j \in \mathbb Z} 
   \Big( \sum _{j' \in \mathbb Z} 2^{-(2-|s|)|j-j'|} 
     \cdot 2^{s j'}\| \phi_{j'} (\Lambda _D)f \|_{L^p } \Big) ^q
   \right\}^{\frac{1}{q}} 
\leq C \| f \|_{\dot B^s_{p,q}} .
\end{split}
\]
Conversely, we have from the resolution by $\{ \psi_{j'} \}_{j' \in \mathbb Z}$ 
that 
\[
\begin{split}
\| f \|_{\dot B^s_{p,q}} 
\leq 
& \left\{ \sum _{ j \in \mathbb Z} 
   \Big( \sum _{j' \in \mathbb Z} 2^{sj} 
     \| \Phi_j(\Lambda _D) \phi_{j} (\Lambda _D) \psi_{j'}(\Lambda _D)f \|_{L^p } \Big) ^q
   \right\}^{\frac{1}{q}} .
\end{split}
\]
For $N \in \mathbb N$ to be fixed later, 
we divide the sume over $j' \in \mathbb Z$ into two cases of 
$|j-j'| \leq N$ and $|j - j'| > N$. 
For the first case, the uniform boudedness of $\phi_j (\Lambda)$ implies that 
a constant $C_N$ depending on $N$ exists such that  
\[
\begin{split}
 \left\{ \sum _{ j \in \mathbb Z} 
   \Big( \sum _{|j-j'| \leq N} 2^{sj} 
     \| \Phi_j(\Lambda _D) \phi_{j} (\Lambda _D) \psi_{j'}(\Lambda _D)f \|_{L^p } \Big) ^q
   \right\}^{\frac{1}{q}}
\leq C_N \| f \|_{\widetilde {\dot B^s_{p,q}}}.
\end{split}
\]
We estimate the second case, 
\[
\begin{split}
& \left\{ \sum _{ j \in \mathbb Z} 
   \Big( \sum _{|j-j'| > N} 2^{sj} 
     \| \Phi_j(\Lambda _D) \phi_{j} (\Lambda _D) \psi_{j'}(\Lambda _D)f \|_{L^p } \Big) ^q
   \right\}^{\frac{1}{q}}
\\
\leq & 
\left\{ \sum _{ j \in \mathbb Z} 
   \Big( \sum _{|j-j'| > N}  \frac{2^{-2j'+2j}}{(1+2^{-2j'+2j})^2}
     \cdot 2^{sj} \|  \phi_{j} (\Lambda _D) f \|_{L^p } \Big) ^q
   \right\}^{\frac{1}{q}}
\leq C \| f \|_{\dot B^s_{p,q}} \sum_{|j'| > N} 2^{-2|j'|}.
\end{split}
\]
The three inequalities above yields that 
\[
\left( 1 - C \sum _{|j'| > N} 2^{-2|j'|}
\right) \| f \|_{\dot B^s_{p,q}} 
\leq C _N \| f \|_{\widetilde {\dot B^s_{p,q}}} ,
\]
and we obtain the converse inequality by taking $N$ sufficiently large. 
\hfill $\Box$ 

\section{Boundary value of $(\nabla ^\perp \Lambda _D^{-1} f \cdot \nabla ) g$}
\label{appen_2}

\noindent 
{\bf Proof of Lemma~\ref{lem:0906-4} (i)}. 
Let $f \in \dot B^0_{\infty,1}$. Then we have from 
Lemma~\ref{lem:0906-8} and $\dot B^0_{\infty,1} \hookrightarrow L^\infty$ that 
\[
f = \sum _{ j \in \mathbb Z} \phi_j(\Lambda _D) 
\quad \text{ in  } L^\infty .
\]
We write 
\[
\begin{split}
\phi_j(\Lambda_D)  f(x) 
=
& (-\Delta _D ) ^{-1} \Big( \Lambda _D^2 \phi_j(\Lambda _D) f\Big) (x)
\\
=
&  \int_{\{ |y| < 1 \}} 
\Big( \frac{1}{2\pi} \log |x-y| + \Phi (x,y)\Big) 
 \Big( \Lambda _D^2 \phi_j(\Lambda _D) f\Big) (y) ~dy ,
\end{split}
\]
and see that 
$\big( \Lambda _D^2 \phi_j(\Lambda _D) f\big)$ is in $L^\infty$ 
by the boundedness of the spectral multiplier \eqref{0830-3}. 
Therefore, the integral above can be regarded as a continuous function 
for $|x| \leq 1$ and it is easy to check that 
the integral is zero when $|x| = 1$. 
\hfill $\Box$

\bigskip 

\noindent 
{\bf Proof of Lemma~\ref{lem:0906-4} (ii)}. 
Since $\Phi(x,y) = \Phi(y,x)$, we write 
\begin{equation}\notag 
2\pi (-\Delta _D)^{-1} (x,y) 
= 2\pi \partial _{x_1} \Big( \dfrac{1}{2\pi}\log |x-y| - \Phi(x,y)\Big) 
= \dfrac{x_1 - y_1}{|x-y|^2} - \dfrac{x_1 - \frac{y_1}{|y|^2}}{|x- \frac{y}{|y|^2}|^2} . 
\end{equation}
When $x = (0,1)$, 
\[
\begin{split}
2\pi \partial _{x_1} \Big( \dfrac{1}{2\pi}\log |x-y| - x \Phi(x,y)\Big) 
=
& \dfrac{- y_1}{|x-y|^2} - \dfrac{ - y_1}{ |y|^2 |x- \frac{y}{|y|^2}|^2}\Big|_{x=(0,1)}
\\
=
& \dfrac{- y_1}{|x-y|^2} - \dfrac{ - y_1}{ |x|^2 |y- \frac{x}{|x|^2}|^2} \Big|_{x=(0,1)}
= 0 .
\end{split}
\]
Therefore, when $x = (0,1)$
\[
(\nabla ^{\perp} \Lambda _D^{-1} f_k \cdot \nabla ) g_l
=\Big(\nabla ^{\perp} (-\Delta _D)^{-1} (\Lambda _D f_k) \cdot \nabla \Big) 
 (-\Delta _D)^{-1}( \Lambda _D ^2 g_l) 
 =0.
\]
We also see that the domain and the derivative 
$\big( \nabla ^\perp \Lambda _D^{-1} (\cdot) \cdot \nabla\big) (\cdot) $ 
are invariant under the rotation 
and we see that $(\nabla ^{\perp} \Lambda _D^{-1} f_k \cdot \nabla ) g_l$ 
is zero on the boundary. 
\hfill $\Box$

\vskip10mm

\noindent
{\bf Acknowledgements. }
The author was supported by the Grant-in-Aid for Young Scientists (A) (No.~17H04824)
from JSPS.

\vskip3mm 
%
%
%

\end{document}